\documentclass[10pt]{amsart}

\usepackage{amsmath,amssymb,amsthm,amsrefs,a4wide}
\usepackage{latexsym}
\usepackage{indentfirst}
\usepackage{graphicx}
\usepackage{placeins}
\usepackage{booktabs}
\usepackage{algorithm}
\usepackage{algorithmic}
\usepackage{multirow}

\usepackage{subfigure}
\usepackage{graphicx,xcolor}
\usepackage{diagbox}
\def\mathrlap#1{\text{\hbox to 5pt{$\mathsurround=0pt#1$\hss\hss}}}

\newcommand{\E}{\mathbb{E}}

\theoremstyle{plain}
\newtheorem{theorem}{Theorem}

\theoremstyle{definition}

\theoremstyle{remark}

\newtheoremstyle{cited}%
  {3pt}
  {3pt}
  {\itshape}
  {}
  {\bfseries}
  {.}
  {.5em}
  {\thmname{#1} \thmnumber{#2} \thmnote{\normalfont#3}}

\theoremstyle{cited}
\newtheorem{citedthm}[theorem]{Theorem}

\DeclareMathOperator{\tr}{tr}

\DeclareMathOperator{\argmin}{argmin}

\newcommand{\bbm}{\begin{bmatrix}}
\newcommand{\ebm}{\end{bmatrix}}

\newcommand{\grad}{\nabla}

\begin{document}

\title[Semigroup method for PDEs and eigenvalue problems]{A semigroup method for high dimensional elliptic PDEs and
  eigenvalue problems based on neural networks}

\author{Haoya Li}\thanks{Department of Mathematics, Stanford University, Stanford, CA 94305, USA
  (\texttt{lihaoya@stanford.edu}).}

\author {Lexing Ying}\thanks{Department of Mathematics and ICME, Stanford University, Stanford, CA
  94305, USA (\texttt{lexing@stanford.edu}).}

\thanks{The work of L.Y. is partially supported by the U.S. Department of Energy, Office of Science,
  Office of Advanced Scientific Computing Research, Scientific Discovery through Advanced Computing
  (SciDAC) program and also by the National Science Foundation under award DMS-1818449.
}


\begin{abstract}
In this paper, we propose a semigroup method for solving high-dimensional elliptic partial
differential equations (PDEs) and the associated eigenvalue problems based on neural networks. For the PDE problems, we reformulate the original equations as variational problems with the help of semigroup operators and then solve the variational problems with neural network (NN) parameterization. The main advantages are that no mixed second-order derivative computation is needed during the stochastic gradient descent training and that the boundary conditions are taken into account automatically by the semigroup operator. Unlike popular methods like PINN \cite{raissi2019physics} and Deep Ritz \cite{weinan2018deep} where the Dirichlet boundary condition is enforced solely through penalty functions and thus changes the true solution, the proposed method is able to address the boundary conditions without penalty functions and it gives the correct true solution even when penalty functions are added, thanks to the semigroup operator. For eigenvalue problems, a primal-dual method is proposed, efficiently resolving the constraint with a simple scalar dual variable and resulting in a faster algorithm compared with the BSDE solver \cite{han2020solving} in certain problems such as the eigenvalue problem associated with the linear Schr\"odinger operator. Numerical results are provided to demonstrate the performance of the proposed methods. 
\end{abstract}

\maketitle

\section{Introduction}\label{sec:intro}

Central to the discipline of applied mathematics is the problem of numerically solving partial differential equations (PDEs), among which the high dimensional problems are particularly challenging due to the ``curse of dimensionality'', a phenomenon that the computational complexity of certain algorithms increases exponentially with the dimension. A more challenging problem is the eigenvalue problem, which is closely related with the PDE problem and suffers from the curse of dimensionality as well. In this paper we limit the discussion to second-order linear PDEs and related eigenvalue problems.

In recent years, deep learning methods have experienced great success across a wide range of domains such as image recognition \cites{10.1145/3065386,he2016deep}, natural language processing \cites{graves2013speech, devlin2018bert}, molecular dynamics simulation \cites{zhang2018deep, jia2020pushing}, and protein structure prediction \cite{alquraishi2019alphafold}. One reason behind this success is that neural network models are good approximators for high-dimensional functions that can be trained efficiently in most cases. Leveraging on this property of the neural network, a myriad of data-driven methods have been proposed for solving high-dimensional PDEs and eigenvalue problems, for example, see \cites{weinan2018deep,li2020solving,raissi2019physics,weinan2017deep}. 

In a recent paper \cite{li2020solving}, the semigroup operator of the differential operator is used to rewrite the variational form, which frees the algorithm from calculations of any mixed second-order derivative and automatically handles the boundary conditions. In this paper, we extend this method to more general elliptic PDEs and also derive a primal-dual method for solving the corresponding eigenvalue problems.

\subsection{Background and related work}
For high-dimensional PDEs, Monte-Carlo methods using Feynman-Kac formulas can be applied to obtain
the value of the approximate solution at a given location. However, satisfactory solutions should
provide information of not only the values on a finite number of points, but also of the entire
landscape. 

Monte Carlo methods have also been widely applied to the eigenvalue problems. Among various
approaches, the variational Monte Carlo method (VMC) and the diffusion Monte Carlo method (DMC) are
two most well-known examples that have been thoroughly investigated in the context of quantum
mechanics, see for example \cite{foulkes2001quantum}. The idea of VMC is to parameterize the wave
function and minimize the energy of the system with respect to the parameters, where the energy is
expressed as an expectation with respect to the probability distribution given by the squared
modulus of the wave function, and is numerically computed via the Monte Carlo method. DMC utilizes
the imaginary-time Schr\"odinger equation, whose solution can be represented by a convolution with
respect to the Green's function and can thus be evaluated by Monte Carlo simulations. Since the
imaginary-time Schr\"odinger equation is a linear differential equation, the component of the lowest
energy eigenfunction remains and other components vanish as the time goes to infinity, and the wave
function of the ground state can be obtained. 

For neural network based approaches, the general idea is to approximate the solution with a
neural network, and then train the neural network to minimize a loss built either from a variational
form of the PDE or from a norm of the residue of an equivalent equation of the original PDE.  Mostly
related to the current work, \cites{khoo2019solving,li2020solving} are concerned with the high
dimensional PDEs describing the committor function in the transition path theory, which is a
second-order elliptic equation with a specific kind of Dirichlet boundary conditions. 

In \cite{weinan2017deep}, the backward stochastic differential equation (BSDE) method forms the
equivalent equation using a BSDE, and the residue norm of the equivalent fixed-point equation is
minimized. In some cases the equivalent variational problem has been established, and we only need
to directly apply the neural network parameterization. For example, in \cite{carleo2017solving}, a
neural network with one hidden layer is used to give the trial functions for the VMC method, and the
gradient function needed in the optimization is also evaluated by Monte Carlo method. In
\cite{weinan2018deep}, a ResNet structure is used to parameterize the approximate solution, and the
optimization problem is obtained from the variational formulation of elliptic PDEs, which is then
solved by stochastic gradient descent (SGD) methods. 

In \cite{han2020solving}, the authors extend the BSDE approach to solve the eigenvalue problem with
a second-order elliptic operator. The differential equation is rewritten as a fixed-point equation
with the help of the corresponding semi-group operator as in the BSDE method. By It\^o's formula,
the semi-group operator is represented by a stochastic integral. After that, the numerical solution is obtained by minimizing the loss defining as the $L^2$ norm of the fixed-point equation residue. The $L^2$ constraint on the eigenfunction is implemented by dividing the $L^2$ norm in each batch during training. 

Since both the proposed method and the BSDE method involves semigroups of diffusion processes, we
remark that there are several major difference between the proposed method and the BSDE
method. Since the BSDE method has many variants, here we take the version in \cite{han2020solving}
to avoid ambiguity. Firstly, the BSDE method usually requires two neural networks, one for the
approximate solution and another for its gradient, while our method only needs a single neural
network for the approximate solution. Secondly, the semigroup used in \cite{han2020solving} is not
the one that corresponds to the second-order differential parameter. In particular, it depends on
$\lambda$, while in our method the semigroup does not depend on $\lambda$.  Finally, in the method
of \cite{han2020solving}, the eigenvalue is also a parameter that needs to be optimized, while in
our method, it does not appear as a variable of the optimization problem but can be computed easily
after solving for the eigenfunction. 

\subsection{Contributions and contents}

The two major contributions of our approach are 
\begin{itemize}
\item The semigroup formulation removes the need of calculating any mixed second-order derivatives,
  and it treats the Dirichlet boundary conditions naturally. Compared with the method proposed in \cite{li2020solving}, the method proposed in this paper applies to problems with non-zero right-hand-side term. Compared with other popular neural network based PDE solvers such as the Deep Ritz method \cite{weinan2018deep} and PINN method \cite{raissi2019physics} where the Dirichlet boundary condition is treated solely through additive penalty functions, the proposed method is capable of addressing the Dirichlet boundary condition without penalty functions, and the ground truth remains unchanged even when penalty functions are added. 
\item We extend the semigroup method also to the eigenvalue problem with a primal dual algorithm,
  which is able to efficiently enforce the constraint on the $L^2$ norm of the approximate solution. In the numerical comparison with the BSDE method proposed in \cite{han2020solving}, the proposed method is shown to have better performance in certain aspects. For example, less running time is needed to achieve the same precision in the linear Sch\"odinger problem. 
\end{itemize}
The rest of the paper is organized as follows. Section \ref{sec:elliptic} describes the semigroup
approach for the second order elliptic PDEs. Section \ref{sec:eigen} discusses the primal-dual
approach for the eigenvalue problems. Finally, numerical results are reported in Section
\ref{sec:numerical}.

\section{Elliptic PDEs}
\label{sec:elliptic}

Consider the following second order elliptic equation:
\begin{equation}\label{eq:elliptic}
-\grad \cdot (a(x)\grad u(x)) = f(x), \quad x\in\Omega
\end{equation}
where the coefficient $a(x)$ is uniformly bounded above zero. Here, we consider two types of
boundary conditions: the Dirichlet boundary condition: 
\begin{equation}\label{eq:dbc}
u(x) = r(x), \quad x\in\partial\Omega,
\end{equation}
and the periodic boundary condition: 
\begin{equation}\label{eq:pbc}
u(x) = u(x+e_i), \quad x\in\Omega = [0,1)^d,
\end{equation}
where $e_i$ is the $i$-th standard basis vector in $\mathbb{R}^d$. When the periodic boundary condition is used, we assume that $a(x)$ and $f(x)$ are both periodic, and we add a further constraint that $\int_{\mathrlap{\Omega}} u(x) \mathrm{d} x = 0$, so that the solution is unique. Otherwise for any solution $u$, $u+C$ is also a solution for an arbitrary constant $C$.

\subsection{Semigroup formulation}
We define $V(x)= -\log(a(x))$ and
let $X_t$ be the solution to the stochastic differential equation (SDE)
\begin{equation}\label{eq:trajectory}
\mathrm{d}X_t = -\grad V(X_t)\mathrm{d}t + \sqrt{2}\mathrm{d}W_t, \quad X_0 = x,
\end{equation}
where $W_t$ is the standard $d$-dimensional Brownian motion. 
For a fixed small time step $\delta>0$, we define the operator $P$ as follows:
\begin{equation}\label{eq:beforesplit}
  (Pu)(x):=\E^{x}\left(u\left(X_{\tau\wedge \delta}\right)\right),
\end{equation}
where $\E^{x}$ is the expectation taken with respect to the law of the process
\eqref{eq:trajectory}, and $\tau = \infty$ if the periodic boundary condition is used, and if the
Dirichlet boundary condition is used, $\tau = \tau_{\partial \Omega}$ is defined as the hitting time
of $\partial\Omega$.  By Dynkin's formula, for the solution $u$ of the equation \eqref{eq:elliptic}
we have
\begin{equation}\label{eq:replace}
    Pu(x)=u(x)+\E^{x}\int_{0}^{\tau\wedge \delta} \mathcal{A} u\left(X_s\right) \mathrm{d} s = u(x)
    - \E^{x}\int_{0}^{\tau\wedge \delta} \frac{f}{a}\left(X_s\right) \mathrm{d} s , \quad \forall x
    \in \Omega,
\end{equation}
where $\mathcal{A} = \Delta -\grad V\cdot \nabla$ is the infinitesimal generator, and thus for the
solution $u$ of the PDE \eqref{eq:elliptic}, we have $\mathcal{A}u(x) = -f(x)/a(x)$. Following
\cite{li2020solving}, $Pu$ can be decomposed into two parts as follows:
\begin{equation}\label{eq:split}
    (P u)(x) = \E^{x}\left(u\left(X_{\tau\wedge \delta}\right)\right) =
  \E^{x}\left(u\left(X_{\delta}\right)\mathbf{1}_{\{\delta<\tau\}}\right) +
  \E^{x}\left(r\left(X_{\tau}\right)\mathbf{1}_{\{\delta\geq \tau\}}\right),
\end{equation}
We denote the first part of \eqref{eq:split} as
\begin{equation}\label{eq:pi}
  (P^i u)(x) \equiv \E^{x}\left(u(X_{\tau\wedge\delta})\mathbf{1}_{\{\delta<\tau\}}\right) =
  \E^{x}\left(u(X_{\delta})\mathbf{1}_{\{\delta<\tau\}}\right),
\end{equation}
where the superscript $i$ stands for the interior contribution and the second part of \eqref{eq:split} as
\begin{equation}
  (P^b r)(x) \equiv
  \E^{x}\left(r(X_{\tau\wedge\delta})\mathbf{1}_{\{\delta\geq\tau\}}\right)=
  \E^{x}\left(r(X_{\tau})\mathbf{1}_{\{\delta\geq\tau\}}\right),
\end{equation}
where the superscript $b$ stands for the boundary contribution. As mentioned earlier, when the periodic boundary condition is used, $\tau = \infty$, so $\mathbf{1}_{\{\delta\geq \tau\}}=0$, and thus $P^b$ becomes a zero operator. In this case, $P^br$ is a zero function for any function $r$. In order to give a uniform formulation for both types of boundary conditions, we set $r(x)=0$ when the periodic boundary condition is used, which has no effect other than making $P^br$ well defined, and any other function can be used. With these operators,
\eqref{eq:replace} can be rewritten succinctly as
\begin{equation}\label{eq:nobd}
  (I-P^i)u(x)-(P^b r)(x) - (Tf)(x) = 0, 
\end{equation}
where $(Tf)(x) = \E^{x}\int_{0}^{\tau\wedge \delta} \frac{f}{a}\left(X_s\right) \mathrm{d} s$. 
This equation can be reformulated as the following variational problem
\begin{equation}\label{eq:variational}
  \begin{aligned}
    \min_u \frac{1}{2}\int_{\mathrlap{\Omega}} u(x)\left((I-P^i)u(x)\right)
    \rho(x) \mathrm{d} x - \int_{\mathrlap{\Omega}} u(x)(P^b r(x)+Tf(x))\rho(x) \mathrm{d} x,
  \end{aligned}
\end{equation}
where $\rho(x) = a(x)/\left(\int_\Omega a(x)\mathrm{d}x\right)$. In order to do this, we need a
result from \cite{li2020solving}.
\begin{citedthm}[\cite{li2020solving}]\label{thm:sym}
  $P^i$ is a symmetric operator on $L_{\rho}^2(\Omega)$, in other words, $\langle u, P^i
  v\rangle_{\rho} = \langle P^i u, v\rangle_{\rho}$, where $\langle \cdot, \cdot\rangle_{\rho}$
  denotes the inner product of the Hilbert space $L_{\rho}^2(\Omega)$.
\end{citedthm}
One can show using Theorem~\ref{thm:sym} that the solution to \eqref{eq:variational} is the same as
the solution to \eqref{eq:nobd} in the following way. Assume that $u^*$ is the solution to
\eqref{eq:variational} and $\eta$ is continuous in $\Omega$ with compact support. By plugging $u(x,
\epsilon) = u^*(x) + \epsilon\eta(x)$ into \eqref{eq:variational} and taking derivative with respect
to $\epsilon$, we obtain
\begin{equation}
  \begin{aligned}
    0 &= \frac{\partial}{\partial \epsilon} \left.\left( \frac{1}{2}\int_{\mathrlap{\Omega}}u(x, \epsilon)\left((I-P^i)u(x, \epsilon)\right) \rho(x) \mathrm{d} x - \int_{\mathrlap{\Omega}} u(x, \epsilon)(P^b r(x)+Tf(x))\rho(x) \mathrm{d} x\right)\right|_{\epsilon = 0}\\
    &=\int_{\mathrlap{\Omega}} \eta(x)\left((I-P^i)u^*(x)\right)\rho(x)\mathrm{d} x
    -\int_{\mathrlap{\Omega}} \eta(x)(P^br(x)+Tf(x))\rho(x) \mathrm{d} x\\
    &=\int_{\mathrlap{\Omega}}\eta(x)\left((I-P^i)u^*(x)-P^b r(x)-Tf(x)\right)\rho(x) \mathrm{d} x,
  \end{aligned}
\end{equation}
and thus $(I-P^i)u^*(x)-P^b r(x)-Tf(x)=0$ since this is true for any $\eta$ that is continuous in
$\Omega$ with compact support.

\subsection{Neural network approximation}
In order to address the curse of dimensionality, the function $u$ in \eqref{eq:variational} is
parameterized with a neural network $u_\theta$, and the optimization problem becomes
\begin{equation}\label{eq:nnvariational}
  \begin{aligned}
    \min_\theta \frac{1}{2}\int_{\mathrlap{\Omega}} u_\theta(x)\left((I-P^i)u_\theta(x)\right)
    \rho(x) \mathrm{d} x - \int_{\mathrlap{\Omega}} u_\theta(x)(P^b r(x)+Tf(x))\rho(x) \mathrm{d} x.
  \end{aligned}
\end{equation}
When the Dirichlet boundary condition \eqref{eq:dbc} is used, we solve the following penalized problem to better address the boundary condition, which is not necessary but is shown to be able to improve the performance of the algorithm. 
\begin{equation}\label{eq:nnvariationalpen}
  \begin{aligned}
    \min_\theta &\int_{\mathrlap{\Omega}} u_\theta(x)\left(\frac{1}{2}(I-P^i)u_\theta(x)-P^b r(x)-Tf(x)\right)\rho(x) \mathrm{d} x + c\int_{\mathrlap{\partial\Omega}} (u_\theta(x)-r(x))^2\mathrm{d}\mu(x),
  \end{aligned}
\end{equation}
where $\mu(x)$ is a probability measure supported on $\partial\Omega$. 

\begin{figure}[H]
  \centering
  \includegraphics[width=0.8\textwidth]{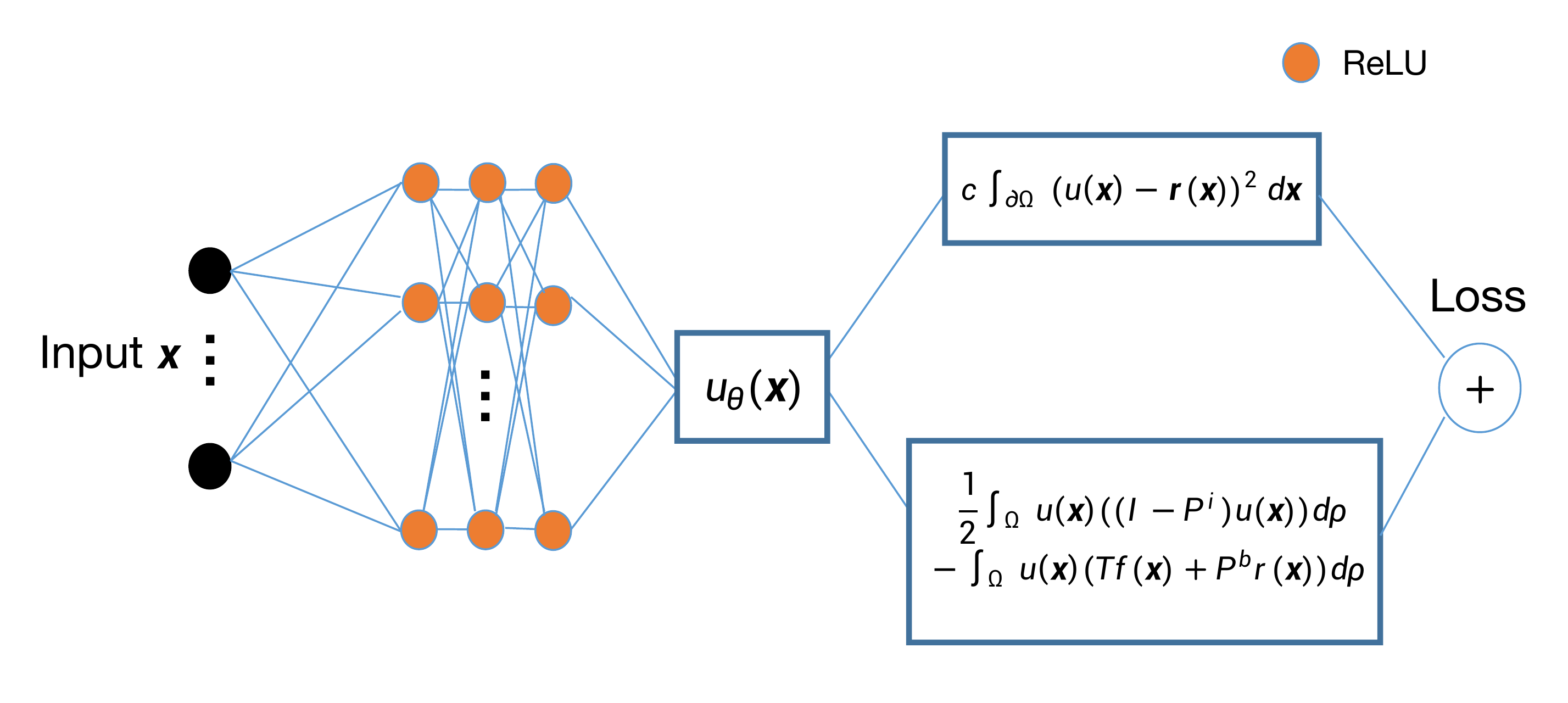}
  \caption{An example of the neural network architecture and the
    corresponding loss when using the Dirichlet boundary condition.
  }
  \label{fig:nn}
\end{figure}

The architecture of the neural network used in this situation is depicted in Figure~\ref{fig:nn}. We
adopt a three layer fully connected neural network with ReLU activation. Compared with
\cite{li2020solving}, we do not have the singularity layer here since we do not have a temperature parameter and the singularities that appear in the situation of extremely high and extermely low temperatures no longer exist.

When the periodic boundary condition \eqref{eq:pbc} is used, we adopt the following neural network
architecture to address the boundary condition with the help of the trigonometric basis, where $m$
is a hyperparameter of the neural network. Since the boundary condition is already treated by the
neural network architecture, we do not need the penalty term in this situation.
\begin{figure}[H]
  \centering
  \includegraphics[width=0.8\textwidth]{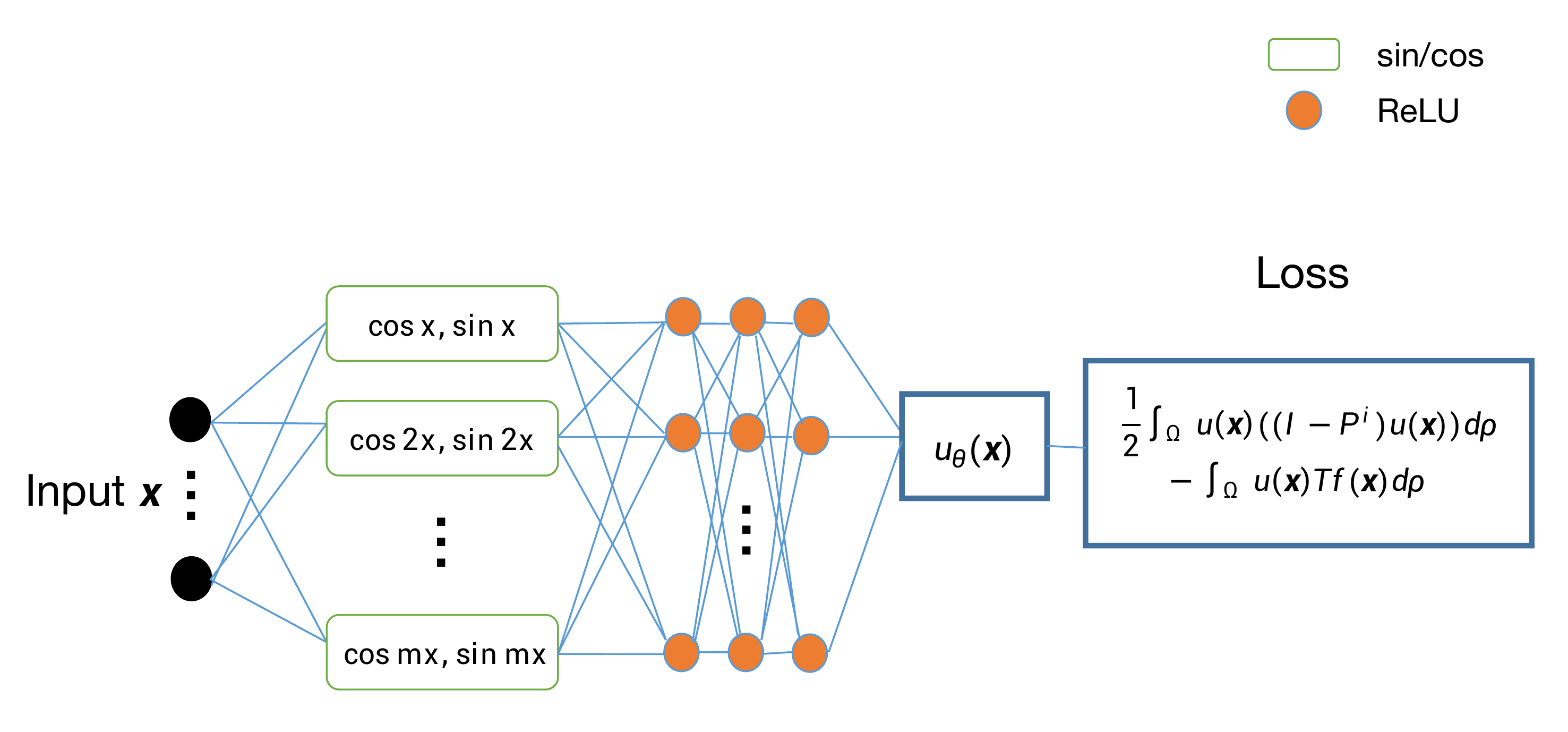}
  \caption{An example of the neural network architecture and the
    corresponding loss when using the periodic boundary condition.
  }
  \label{fig:nnperiod}
\end{figure}

In the implementation of the optimization algorithm, we also need the derivative of the integrals in
\eqref{eq:nnvariational}. By the symmetry of $P^i$ in $L_{\rho}^2(\Omega)$, the derivative is
\begin{equation}
  \int_{\mathrlap{\Omega}} \grad_\theta u_\theta(x)\left((I-P^i)u_\theta(x)\right)
  \rho(x) \mathrm{d} x - \int_{\mathrlap{\Omega}} \grad_\theta u_\theta(x)(P^b r(x)+Tf(x))\rho(x) \mathrm{d} x.
\end{equation}
With the help of a random variable $X\sim \rho$, the derivative can be further transformed into:
\begin{equation}\label{eq:gradexp}
  \E_{X\sim \rho} \grad_\theta u_\theta(X)\left((I-P^i)u_\theta(X)- (P^b r(X)+Tf(X))\right).
\end{equation}
An unbiased estimator for \eqref{eq:gradexp} is thus
\begin{equation}\label{eq:unbiased1}
    \grad_\theta u_\theta(X)\left(u_\theta(X)-u_\theta(X_\delta)\mathbf{1}_{\{\delta<\tau\}}-
    r(X_\tau)\mathbf{1}_{\{\delta\geq\tau\}}-\int_{0}^{\tau\wedge \delta}
    \frac{f}{a}\left(X_s\right) \mathrm{d} s\right),
\end{equation}
where $X_s$ is the solution of \eqref{eq:trajectory} at time $s$ with the initial condition given by
$X_0=X$. $X_\delta$ and $X_\tau$ are obtained by evaluating $X_s$ at $\delta$ and the stopping time
$\tau$, respectively. The derivative of the penalty term is
\begin{equation}\label{eq:unbiasedbd}
    2c\int_{\mathrlap{\partial\Omega}} \grad_\theta u_\theta(x)(u_\theta(x)-r(x))\mathrm{d}\mu(x)
\end{equation}
and an unbiased estimator for the derivative is
\begin{equation}\label{eq:unbiasedbd1}
    2c\grad u_\theta(X)(u_\theta(X)-r(X)),
\end{equation}
where $X\sim \mu$. Notice that there is no mixed second-order derivative in \eqref{eq:unbiased1} and \eqref{eq:unbiasedbd1}, unlike \cite{khoo2019solving} or \cite{weinan2018deep}, where mixed second-order derivatives $\frac{\partial^2 }{\partial \theta\partial x_i}$, $1\leq i\leq d$ are needed in the training process. 

The optimization problem \eqref{eq:nnvariational} and the penalized problem
\eqref{eq:nnvariationalpen} can be solved by applying SGD-type optimization, for example the Adam
method in \cite{kingma2014adam}. In the numerical implementation, the integral $\int_{0}^{\tau\wedge
  \delta} \frac{f}{a}\left(X_s\right) \mathrm{d} s$ can be approximated for example by the
Euler-Maruyama scheme (see for example \cite{kloeden2013numerical}). For other implementation
details such as the determination of $\mathbf{1}_{\{\delta\geq\tau\}}$, we follow the method in
\cite{li2020solving}. The complete algorithm is summarized in Algorithm~\ref{alg:semi1}.

\begin{algorithm}[htbp]
    \caption{Semigroup method for the elliptic PDE\eqref{eq:elliptic}}
    \label{alg:semi1}
    \begin{algorithmic}[1]
        \REQUIRE batch size $B$ and $\tilde{B}$, total number of iterations $T_{\mathrm{iter}}$,
        time step $\delta$, learning rate $\eta_t$, the penalty coefficient $c$ used in the Dirichlet case, the level $m$ of trigonometric bases used in the periodic case.\\ 
        \STATE Initialize the neural network $u_\theta$.

        \FOR {$t=1,..,T_{\mathrm{iter}}$}

        \STATE Sample a batch of data from the distribution $\rho$ with size $B$:
        \vspace{-0.4em}
        $$X_{1}, X_{2}, \ldots, X_{B} \sim \rho.$$
        \vspace{-1.5em}

        \STATE For each $X_k$ ($1\leq k\leq B$), sample $X_{k,\delta}$ according to the SDE
        \eqref{eq:trajectory}:
        \vspace{-0.5em}
        $$X_{k, \delta} = X_k - \grad V(X_k)\delta + \sqrt{2}W_{\delta},$$ and when the periodic
        boundary condition is used, move $X_{k,\delta}$ into $\Omega$ by translating an integer
        multiple of the period.
        
        \STATE For each $X_k$, decide the value of $\mathbf{1}_{\{\delta<\tau\}}$ by
        $$\mathbf{1}_{\{\delta<\tau\}} = 1 ~\text{if}~ X_{k,\delta}\in\Omega.$$

        \STATE For each $X_k$ such that $\mathbf{1}_{\{\delta<\tau\}}=0$, let $X_\tau$ be the
        intersection of $\partial \Omega$ and the line segment $X_kX_{k,\delta}$.

        \STATE Compute the gradient in \eqref{eq:unbiased1}. 

        \STATE If the Dirichlet boundary condition is used, sample a batch of data
        $\{\tilde{X}_j\}_{j=1}^{\tilde{B}}$ from the distribution $\mu$, and compute the gradient of
        the penalty term by \eqref{eq:unbiasedbd}.

        \STATE Update the neural network parameters $\theta$ via the Adam method with learning rate $\eta_t$
        (other hyper-parameters in the Adam method are set as the default values).

        \ENDFOR
        \STATE If the periodic boundary condition is used, then sample a batch of data
        $\{\tilde{X}_j\}_{j=1}^{\tilde{B}}$ from the uniform distribution on $\Omega$, and subtract $\frac{1}{\tilde{B}}\sum_{j=1}^{\tilde{B}}u_\theta(\tilde{X}_j)$ from $u_\theta$. 
    \end{algorithmic}
    \vspace{0.5em}
\end{algorithm}

\section{Primal-dual formulation for eigenvalue problems}
\label{sec:eigen}
In this section we extend the semigroup method described in Section \ref{sec:elliptic} to a
primal-dual algorithm for eigenvalue problems. For the eigenvalue problems, a major difference is
that the corresponding variational problem has a constraint $\|u\| = 1$. In \cite{han2020solving},
the authors divide the neural network approximate solution by a normalization factor when
calculating the training loss, which is similar with the Batch Normalization technique, except that
the normalization factor is computed with an auxiliary batch rather than the original batch. Here we
propose a primal-dual method that handles the constraint via a scalar Lagrange multiplier. As a
result, the time complexity of solving the eigenvalue problem is only marginally higher than that of
solving the corresponding PDE.

Consider $\Omega = [0, 1)^d$ and define $C_{\text{per}}^\infty$ as the space of smooth functions on $\Omega$ satisfying the the periodic boundary condition \eqref{eq:pbc}. Consider a symmetric elliptic operator $L = -\sum_{i,j=1}^d\frac{\partial}{\partial x_j}\left(A_{ij}\frac{\partial}{\partial x_i}\right) + V$, where $V\in C_{\text{per}}^{\infty}(\Omega)$, $A_{ij}\in C_{\text{per}}^{\infty}(\Omega),~ 1\leq i, j \leq d$ and $A_{ij} = A_{ji}, ~1\leq i, j \leq d$. We also assume that $A$ is uniformly elliptic, i.e., $\sum_{i,j=1}^d A_{ij}(x)\xi_i\xi_j\geq \alpha \|\xi\|^2$ for some constant $\alpha>0$ and any $x, \xi\in \Omega$. By variational principle we know that for the first eigenpair $(u, \lambda)$, the eigenvalue problem 
\begin{equation}
  Lu = \lambda u, \quad \|u\|_{L^2} = 1,
\end{equation}
with the periodic boundary condition \eqref{eq:pbc} is equivalent to the variational form:
\begin{equation}\label{eq:eigenvargen}
\begin{aligned}
  u &= \underset{u\in H_{\text{per}}^1(\Omega),\|u\|_{L^2}=1}{\argmin} ~\frac{1}{2}\int_{\mathrlap{\Omega}} u(x)Lu(x) \mathrm{d} x,
\end{aligned}
\end{equation}
where $H_{\text{per}}^1(\Omega)$ consists of the functions satisfying the periodic boundary condition \eqref{eq:pbc} in $H^1(\Omega)$. For simplicity of the notations, we omit the function space when there is no ambiguity. 

\subsection{Resolving the constraint}
The main obstacle in implementing the formulation \eqref{eq:eigenvargen} lies in the constraint $\|u\|=1$. We propose to handle the constraint with a multiplier term, and reconstruct \eqref{eq:eigenvargen} as the minimax formulation:

\begin{equation}\label{eq:eigendual}
  \underset{u}{\min}~\underset{g}{\max} ~\frac{1}{2}\int_{\mathrlap{\Omega}} u(x)Lu(x) \mathrm{d} x + \frac{g}{2}(\|u\|_{L^2}^2-1),
\end{equation}
where $g$ is the Lagrange multiplier. The equivalence between \eqref{eq:eigenvargen} and
\eqref{eq:eigendual} can be seen by maximizing over $g$ in \eqref{eq:eigendual} and thus we have
removed the explicit constraint in \eqref{eq:eigenvargen}. We adopt the neural network
parameterization as in Section \ref{sec:elliptic} and obtain the following optimization problem:
\begin{equation}\label{eq:nndual}
  \underset{\theta}{\min}~\underset{g}{\max} ~\frac{1}{2}\int_{\mathrlap{\Omega}}
  u_\theta(x)Lu_\theta(x) \mathrm{d} x + \frac{g}{2}(\|u_\theta\|_{L^2}^2-1).
\end{equation}
Let us define
\begin{equation}
E(\theta,g) = \frac{1}{2}\int_{\mathrlap{\Omega}} u_\theta(x)Lu_\theta(x) \mathrm{d} x +
\frac{g}{2}(\|u_\theta\|_{L^2}^2-1).
\end{equation}
Taking derivatives with respect to $\theta$ and $g$ gives
\begin{equation}
  \begin{aligned}
    \frac{\partial E}{\partial \theta}& = \int_{\mathrlap{\Omega}} \grad_\theta u_\theta(x) (L+g)u_\theta(x) \mathrm{d} x,\\
    \frac{\partial E}{\partial g}& = \frac{1}{2}(\|u_\theta\|_{L^2}^2-1),
  \end{aligned}
\end{equation}
where we have utilized the symmetry of $L$. 

\subsection{Scaling the Lagrange multiplier}
In the optimization problem \eqref{eq:nndual}, it is important that the Lagrange multiplier has an
appropriate scale, since otherwise the problem would be ill-conditioned. If we introduce a scaling
parameter $c$ for the dual variable $g$, then the optimization problem becomes
\begin{equation}\label{eq:minimaxscale}
  \underset{\theta}{\min}~\underset{g}{\max} ~\frac{1}{2}\int_{\mathrlap{\Omega}}
  u_\theta(x)Lu_\theta(x) \mathrm{d} x + \frac{cg}{2}(\|u_\theta\|_{L^2}^2-1),
\end{equation}
or equivalently, 
\begin{equation}
  E = \frac{1}{2}\int_{\mathrlap{\Omega}} u_\theta(x)Lu_\theta(x) \mathrm{d} x + \frac{cg}{2}(\|u_\theta\|_{L^2}^2-1),
\end{equation}
and the derivatives are scaled accordingly:
\begin{equation}\label{eq:dualgrad}
  \begin{aligned}
    \frac{\partial E}{\partial \theta}& = \int_{\mathrlap{\Omega}} \grad_\theta u_\theta(x) (L+cg)u_\theta(x) \mathrm{d} x,\\
    \frac{\partial E}{\partial g}& = \frac{c}{2}(\|u_\theta\|_{L^2}^2-1). 
  \end{aligned}
\end{equation}
If $c$ is too large, then in the optimization process the first term in $E$ is neglected by the
neural network. If $c$ is too small, then the constraint is not well-enforced. Therefore in the
implementation, $c$ should be properly chosen such that the two terms in $E$ are balanced.

Based on \eqref{eq:dualgrad}, we arrive at the primal-dual scheme
\begin{equation}\label{eq:gradscheme}
\begin{aligned}
\dot{\theta}& = -\int_{\mathrlap{\Omega}} \grad_\theta u_\theta(x) (L+cg)u_\theta(x)\mathrm{d} x,\\
\dot{g}& = \frac{c}{2}(\|u_\theta\|_{L^2}^2-1),
\end{aligned}
\end{equation}
In what follows, we focus on the Schr\"odinger operator $L = -\Delta + V$ and consider the problem
of finding the first eigenpair (which corresponds to the ground state) with periodic boundary
conditions. The minimax formulation \eqref{eq:minimaxscale} becomes
\begin{equation}\label{eq:schroddual}
  \underset{u}{\min}~\underset{g}{\max} ~\frac{1}{2}
  \int_{\mathrlap{\Omega}} |\grad u(x)|^2  \mathrm{d} x + \frac{1}{2}\int_{\mathrlap{\Omega}} V(x) |u(x)|^2 \mathrm{d} x + \frac{cg}{2}(\|u\|_{L^2}^2-1).
\end{equation}
Two stochastic schemes for implementating \eqref{eq:gradscheme} are discussed below.

\subsection{Scheme I}
By replacing the operator $\frac{\Delta}{2}$ with the semigroup approximation, we can reformulate
\eqref{eq:schroddual} as
\begin{equation}\label{eq:eigenprimaldual1}
  \underset{u}{\min}~\underset{g}{\max} ~\frac{1}{\delta}\int_{\mathrlap{\Omega}} u(x)(u(x)-\E(u(x +
  W_\delta)))\mathrm{d} x + \frac{1}{2}\int_{\mathrlap{\Omega}} V(x) |u(x)|^2 \mathrm{d} x +
  \frac{cg}{2}(\|u\|^2-1),
\end{equation}
where $W_\delta$ is the standard Brownian motion. This replacement can be justified as follows. By Dynkin's formula, we have
\[
\frac{1}{\delta}(u(x)-\E u(x+W_\delta)) = \E\frac{1}{\delta}\int_0^\delta -\frac{1}{2}\Delta u(x+W_s)\mathrm{d}s.
\]
Let $\varphi = -\frac{1}{2}\Delta u$. If $\varphi$ is bounded and sufficiently smooth, then by Fubini's theorem,
\[
\begin{aligned}
\E\frac{1}{\delta}\int_0^\delta \varphi(x+W_s)\mathrm{d}s &= \frac{1}{\delta}\int_0^\delta \E\varphi(x+W_s)\mathrm{d}s = \frac{1}{\delta}\int_0^\delta \E\Big[\varphi(x) + \nabla\varphi(x)\cdot W_s + O(\|W_s\|^2)\Big]\mathrm{d}s \\
&= \frac{1}{\delta}\int_0^\delta \Big[\varphi(x) + O(s)\Big]\mathrm{d}s = \varphi(x) + O(\delta),
\end{aligned}
\]
which shows that $\frac{1}{\delta}(u(x)-\E u(x+W_\delta)) = -\frac{1}{2}\Delta u(x) + O(\delta)$ and justifies the replacement in \eqref{eq:eigenprimaldual1}. 
Consider a uniform random variable $X$ on
$\Omega$, the problem \eqref{eq:eigenprimaldual1} can be written as
\begin{equation}\label{eq:eigenprimalduale}
  \underset{u}{\min}~\underset{g}{\max} ~\frac{1}{\delta}~\E u(X)(u(X)-u(X + W_\delta)) +
  \frac{1}{2}\E V(X) |u(X)|^2 + \frac{cg}{2}(\|u\|^2-1).
\end{equation}
Since $X$ is uniform on $\Omega$, the process $X+W_\delta$ is reversible and
\begin{equation}
    \langle u,  \tilde{P} v\rangle = \langle \tilde{P} u, v\rangle,
\end{equation}
where $\langle u, v\rangle = \E u(X)v(X)$ and $\tilde{P}u(x) = \E u(x+W_\delta)$. By this symmetry,
a gradient descent scheme similar with \eqref{eq:gradscheme} can be derived
\begin{equation}\label{eq:gradschemesemi}
\begin{aligned}
\dot{\theta}& = -\E\grad_\theta u_\theta(X) \left(\frac{1}{\delta}(u_\theta(X)-u_\theta(X+W_\delta))+cgu_\theta(X)\right),\\
\dot{g}& = \frac{c}{2}\E (u_\theta(X)^2-1). 
\end{aligned}
\end{equation}
Unbiased estimators for the expectations in \eqref{eq:gradschemesemi} are therefore
\begin{equation}\label{eq:estimator1}
  \grad_\theta u_\theta(X) (\frac{1}{\delta}(u_\theta(X)-u_\theta(X+W_\delta))+cgu_\theta(X)),
\end{equation}
and 
\begin{equation}\label{eq:estimatordual}
  \frac{c}{2}(u_\theta(X)^2-1),
\end{equation}
where $X$ is uniform in $\Omega$. 

\subsection{Scheme II}
In \eqref{eq:eigenprimalduale}, we have used $\frac{1}{\delta}(u(x)-\E u(x+W_\delta))$ to approximate
$-\frac{\Delta}{2}u(x)$. The effectiveness of this approximation relies on the spatial symmetry of
$W_\delta$, which is only true in distribution. In other words, after using the unbiased estimators
given in \eqref{eq:estimator1} with samples of $X$ and $W_\delta$, the symmetry can be
impacted. This problem can be addressed by starting by approximating the first term of
\eqref{eq:schroddual} with

\begin{equation}\label{eq:scheme2first}
  \begin{aligned}
    \frac{1}{2}\int_{\mathrlap{\Omega}} |\grad u(x)|^2 \mathrm{d} x &=
    \frac{1}{2\delta}\left(\int_{\mathrlap{\Omega}} \E|u(x)-u(x+W_\delta)|^2 \mathrm{d} x +
    o(\delta)\right).
  \end{aligned}
\end{equation}
The formulation above can be obtained, for example, by plugging $u(x+W_\delta) = u(x)+\grad u(x)\cdot W_\delta+O(\|W_\delta\|^2)$ into the expectation, which gives
\begin{equation}
\begin{aligned}
    &\E|u(x)-u(x+W_\delta)|^2 = \E|\grad u(x)\cdot W_\delta|^2 + o(\delta)= \E (W_\delta^\top\grad u(x)(\grad u(x))^\top W_\delta) + o(\delta)\\
    &= \E(\tr(W_\delta^\top\grad u(x)(\grad u(x))^\top W_\delta)) + o(\delta)= \E(\tr(\grad u(x)(\grad u(x))^\top W_\delta W_\delta^\top)) + o(\delta)\\
    &= \tr(\grad u(x)(\grad u(x))^\top \E(W_\delta W_\delta^\top)) + o(\delta) = \delta\tr(\grad u(x)(\grad u(x))^\top) + o(\delta)\\
    &= \delta|\grad u(x)|^2 + o(\delta). 
\end{aligned}
\end{equation}
Notice that compared with the first term in \eqref{eq:eigenprimaldual1}, the right hand side of \eqref{eq:scheme2first} has an extra factor $1/2$. The appearance of this factor is natural since
\[
\frac{1}{2}\E|u(x)-u(x+W_\delta)|^2 = u(x)(u(x)-\E u(x+W_\delta)) + \frac{1}{2}(\E u(x+W_\delta)^2-u(x)^2),
\]
so the additional factor $1/2$ here can be viewed as a result of replacing $\frac{1}{2}u(x)^2$ with $\frac{1}{2}\E u(x+W_\delta)^2$ in \eqref{eq:eigenprimaldual1}. 

With this approximation and the uniform random variable $X$ on $\Omega$, \eqref{eq:schroddual} can
be written as
\begin{equation}\label{eq:eigenprimaldual2}
  \underset{u}{\min}~\underset{g}{\max}~\E\left(\frac{1}{2\delta}|u(X)-u(X+W_\delta)|^2  + \frac{1}{2} V(X) |u(X)|^2 + \frac{cg}{2}(u(X)^2-1)\right). 
\end{equation}
For this problem, the gradient descent scheme is
\begin{equation}\label{eq:gradschemerand}
\begin{aligned}
  \dot{\theta}& = -\E \bigg[~\frac{1}{\delta}(\grad_\theta u_\theta(X)-\grad_\theta u_\theta(X+W_\delta))(u_\theta(X)-u_\theta(X+W_\delta))\\
    &\quad\quad\quad+ \grad_\theta u_\theta(X)(V(X)u_\theta(X)+cgu_\theta(X))\bigg],\\
  \dot{g}& = \frac{c}{2}~\E (u_\theta(X)^2-1). 
\end{aligned}
\end{equation}
Unbiased estimators for the expectations are therefore
\begin{equation}\label{eq:estimator2}
  \begin{aligned}
    ~\frac{1}{\delta}(\grad_\theta u_\theta(X)-\grad_\theta
    u_\theta(X+W_\delta))(u_\theta(X)-u_\theta(X+W_\delta))+ \grad_\theta
    u_\theta(X)(V(X)u_\theta(X)+cgu_\theta(X))
  \end{aligned}
\end{equation}
and
\begin{equation}
  \frac{c}{2}(u_\theta(X)^2-1).
\end{equation}
It is clear that even after replacing the expectation with the unbiased estimators, the problem is
still symmetric.

\begin{algorithm}[htbp]
    \caption{Semigroup methods for the eigenvalue problem.}
    \label{alg:semi2}
    \begin{algorithmic}[1]
      \REQUIRE default value of the dual variable $g_{\text{default}}$, batch size $B$ and $\tilde{B}$, total number of iterations $T_{\mathrm{iter}}$, time step $\delta$, learning rate $\eta_t, \tilde{\eta}_t$, scaling factor $c$, the level $m$ of trigonometric bases used.\\
      \STATE Initialize the neural network $u_\theta$.
      \FOR {$t=1,..,T_{\mathrm{iter}}$}
      \STATE Sample a batch of data $\{X_{k}\}_{k=1}^B$ uniformly in $\Omega$. 
      \STATE Compute the gradient in \eqref{eq:estimator1} or \eqref{eq:estimator2}. 
      \STATE Sample a batch of data $\{\tilde{X}_j\}_{j=1}^{\tilde{B}}$ uniformly in $\Omega$. 
      \STATE Compute $\epsilon_t = \max\left(\frac{1}{\tilde{B}}\sum_j u_\theta(\tilde{X}_j)^2-1, 1\right)$. 
      \IF{$t=1$ \OR $\epsilon_t\epsilon_{t-1}<0$}
      \STATE Set $g = \mathrm{sign}(\epsilon_t)g_{\text{default}}$. 
      \ELSE 
      \STATE Update the dual variable by $g\leftarrow g + \tilde{\eta}_t\frac{c\epsilon_t}{2}$. 
      \ENDIF
      \STATE Update the neural network parameters by the Adam method with learning rate $\eta_t$ (other hyper-parameters in the Adam method are set as the default values).
      \ENDFOR
    \end{algorithmic}
    \vspace{0.5em}
\end{algorithm}

\subsection{Implementation of the multiplier term}
Intuitively, the dual variable should be able to give the correct preference for
$u_\theta$. Specifically, when $\|u_\theta\|>1$, $g$ should be positive, and when $\|u_\theta\|<1$,
$g$ should be negative. In the implementation, this is enforced in each step of the (stochastic)
gradient update. More specifically, in step $t$ we estimate $\|u_\theta\|-1$ using a batch of data with size $\tilde{B}$, and denote the estimator as $\epsilon_t$. If $\epsilon_t\epsilon_{t-1}<0$, then we reset the multiplier $g$ as $\mathrm{sign}(\epsilon_t)g_{\text{default}}$, where $g_{\text{default}}>0$ is a hyperparameter. In this way, $g$ will be positive when $\|u_\theta\|-1$ changes from negative to positive, and negative when $\|u_\theta\|-1$ changes from positive to negative, which leads $\|u_\theta\|$ towards the correct update direction.

We also set the estimation of the gradient of $g$ to be $c/2$ if it is larger than
$c/2$. This is because $\frac{c}{2}~\E (u_\theta(X)^2-1)$ is lower bounded by $-c/2$ but it has no
upper bound. In practice setting a symmetric upper bound usually gives a better performance. We
summarize the implementation details in Algorithm~\ref{alg:semi2}, where we adopt the neural network
architecture in Figure~\ref{fig:nnperiod}.

\section{Numerical experiments}
\label{sec:numerical}

\subsection{Elliptic PDEs}
In this section, we verify the effectiveness of Algorithm~\ref{alg:semi1}. We measure the error of the numerical solutions by $E_0 = \|u_{\theta}-u^*\|_{L_\rho^2(\Omega)}/\|u^*\|_{L_\rho^2(\Omega)}$, where $u^*$ is the ground truth. In the numerical examples, $E_0$ is estimated on a test set of size $1\times 10^4$.  

\begin{figure}[htbp]

  \centering
  \subfigure[The error of the numerical solution in the process of training. \label{fig:test2_train}]{
    \includegraphics[width=0.4\textwidth]{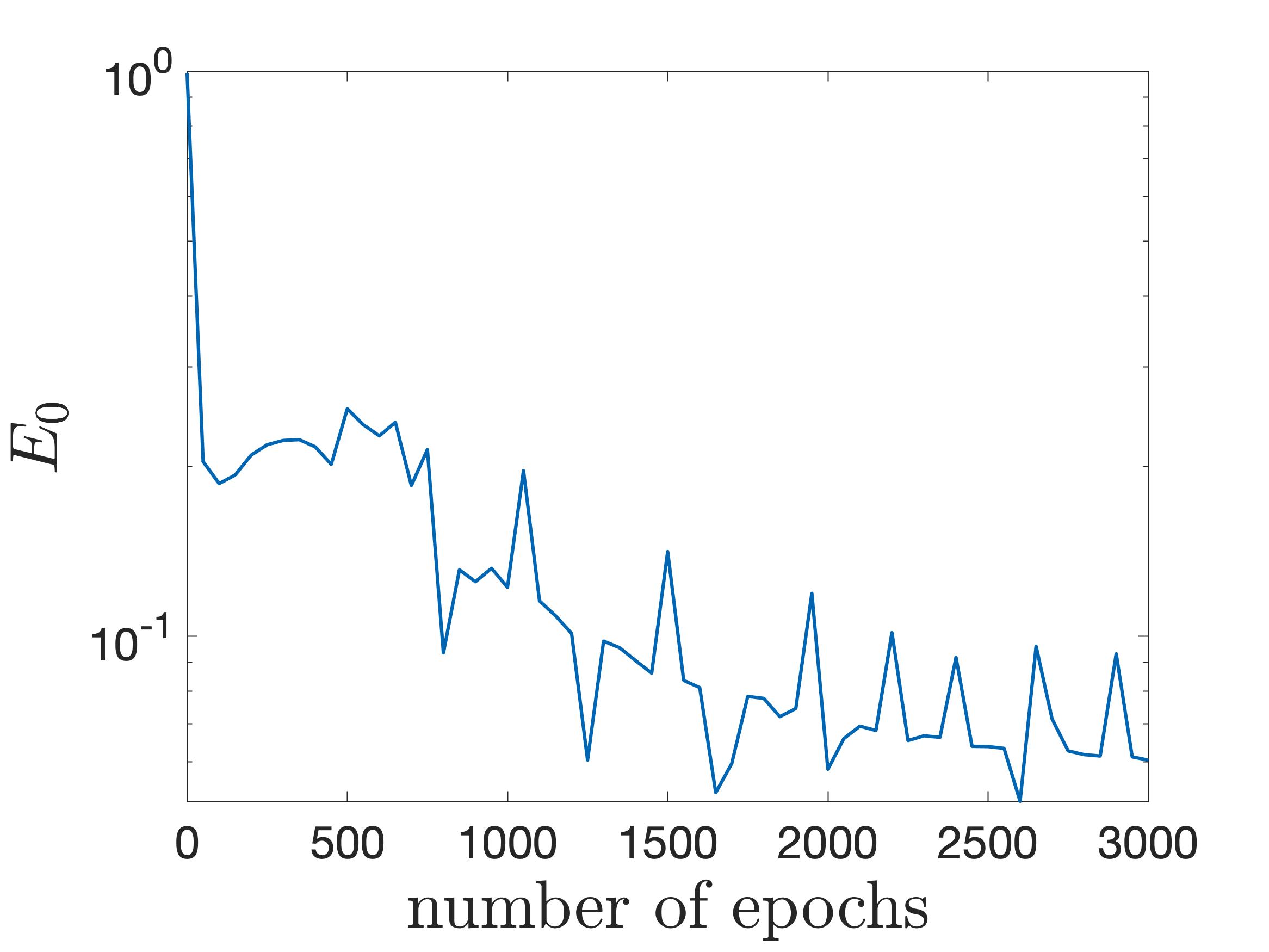} }
  \hspace{0.7em} \subfigure[Comparison of the numerical solution and the ground truth along the $x_4$ axis. \label{fig:radial}]{
    \includegraphics[width=0.4\textwidth]{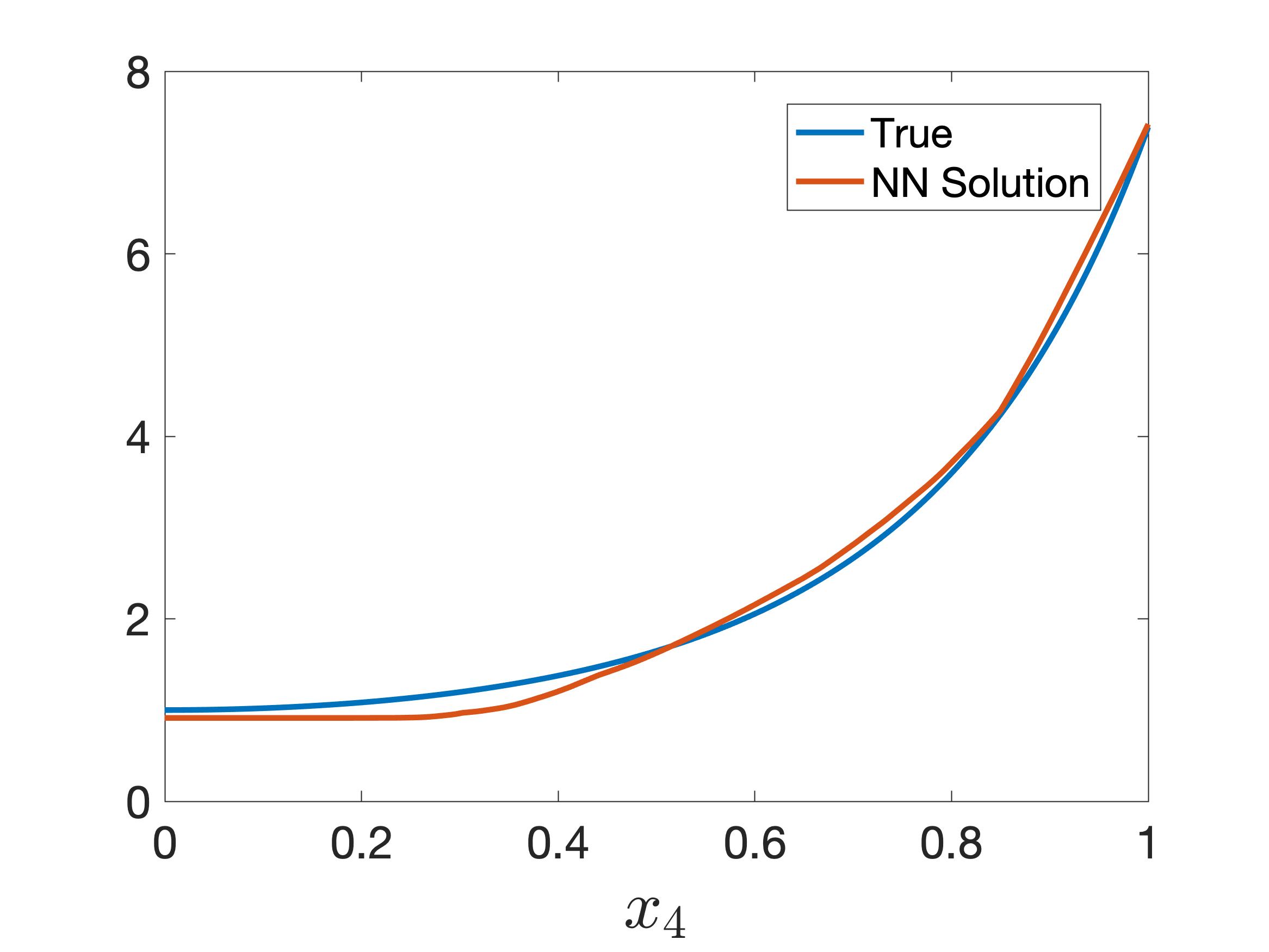} }
  \caption{The convergence of the error $E_0$ and the comparison of the NN solution with the true solution. (a): The error $\|u_{\theta}-u^*\|_{L_\rho^2(\Omega)}/\|u^*\|_{L_\rho^2(\Omega)}$ every $50$ epochs during the training process. (b): Comparison of the NN solution with the true solution on a randomly chosen coordinate axis. Blue: The true solution; Red: The NN solution. \label{fig:test2}}
\end{figure}

In the first numerical example, we set $a(x) = e^{-2\|x\|^2}$ and $f(x) = -4d$ in \eqref{eq:elliptic}, and assume the Dirichlet boundary condition \eqref{eq:dbc} with $r(x) = e^{2}$ on the boundary of the domain $\Omega = B(0,1)$, the unit ball in $\mathbb{R}^d$. In other words, we consider the following problem:
\[
\begin{aligned}
\Delta u + 4x\cdot\nabla u &= -4de^{2\|x\|^2}, \quad x \in B(0, 1),\\
u &= e^2, \quad x\in \partial B(0, 1).
\end{aligned}
\]
In this problem, the unique solution is $u^*(x) = e^{2\|x\|^2}$, which can be obtained by direct calculations. Numerical experiment is implemented for the $10$-dimensional case. We use $4.0\times 10^6$ samples to train the model, which is a three-layer fully-connected network with $120$ neurons in each layer. We set $\delta = 0.0001$, $c=0.8$, and train the model for $3000$ epochs using the Adam algorithm in \cite{kingma2014adam} with learning rate $0.001$ and $70000$ samples in each batch. The final error is $E_0 = 0.060$, which is computed on a test set with size $1.0\times10^5$. The result is visualized in Figure~\ref{fig:test2}. In Figure~\ref{fig:test2_train}, the error is recorded every $50$ epochs, from which we can see that the error converges through the training process. Since the true solution is a radial function, we can visualize the difference between the NN solution and the true solution on a randomly chosen coordinate axis. In Figure~\ref{fig:radial}, it can be seen that the NN solution shows good accordance with the ground truth on the positive $x_4$ axis. From this numerical test, we have verified that the proposed method is effective when the right-hand-side term is non-zero, which is a case not covered in \cite{li2020solving}.

\begin{figure}[htbp]
  \centering
  \subfigure[The error of the numerical solution in the process of training. \label{fig:testper_train}]{
    \includegraphics[width=0.4\textwidth]{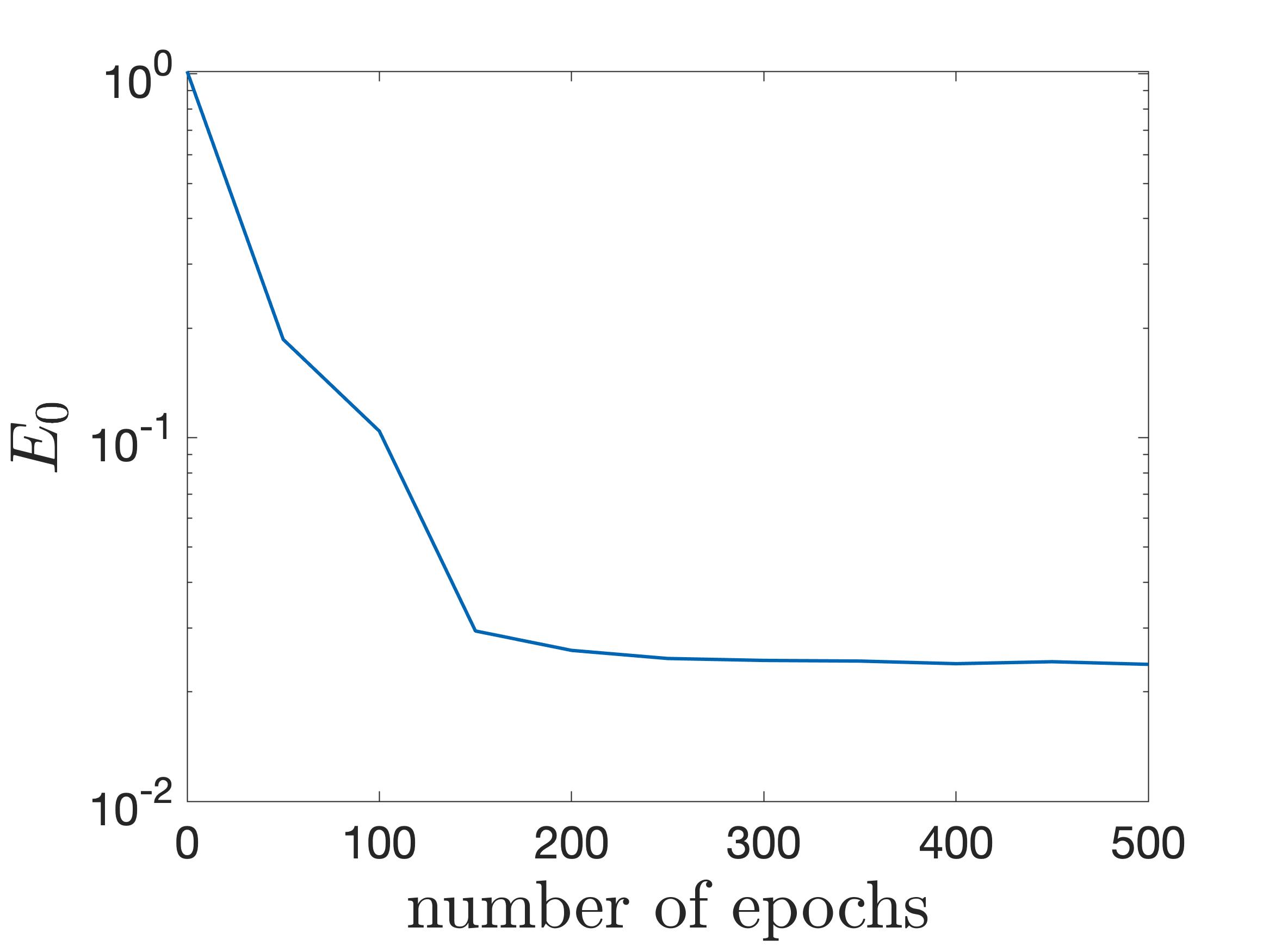} }
  \hspace{0.7em} \subfigure[Comparison of the numerical solution and the ground truth along the $x_4$ axis. \label{fig:randaxis}]{
    \includegraphics[width=0.4\textwidth]{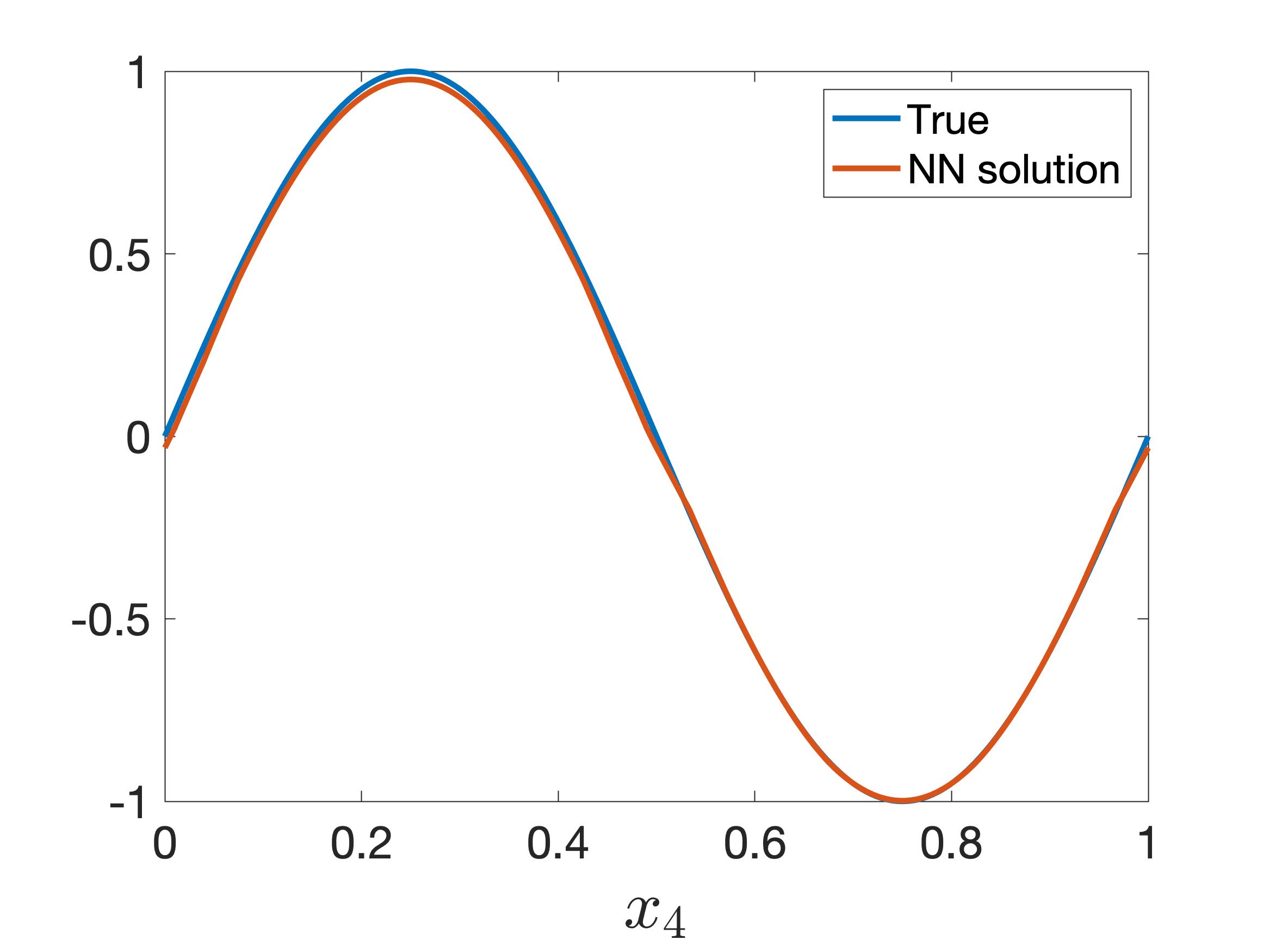} }
  \caption{The convergence of the error $E_0$ and the comparison of the NN solution with the true solution. (a): The error $\|u_{\theta}-u^*\|_{L_\rho^2(\Omega)}/\|u^*\|_{L_\rho^2(\Omega)}$ every $50$ epochs during the training process. (b): Comparison of the NN solution with the true solution on a randomly chosen coordinate axis. Blue: The true solution; Red: The NN solution. \label{fig:testper}}
\end{figure}

In the second numerical example, we consider a problem with the periodic boundary condition \eqref{eq:pbc}. The
parameters are
\begin{equation}
  \begin{aligned}
    a(x) &= \exp\left(-\sum_{i=1}^d \cos(2\pi x_i)\right),\\
    f(x) &= 2\pi^2 \exp\left(-\sum_{i=1}^d \cos(2\pi x_i)\right)\left(\sum_{i=1}^d \left(2\sin(2\pi x_i) -  \sin(4\pi x_i) \right)\right).
  \end{aligned}
\end{equation}
By direct calculation, we see that the exact solution is
\begin{equation}
\begin{aligned}
  u(x) &= \sum_{i=1}^d \sin(2\pi x_i).
\end{aligned}
\end{equation}
Numerical experiment is carried out for $d = 10$. The neural network structure is descibed in Figure~\ref{fig:nnperiod}, in which we take $m=1$ and set the width of the network as $12$. We take $\delta=0.0001$, $\tilde{B} = 2\times 10^5$. In the implementation of the Adam algorithm, the batchsize is set as $B=70000$ and the learning rate is set as $\eta=0.001$. A training set with $1.0\times 10^7$ samples to train the model. After a training
process of $500$ epochs, the final error
$\|u_{\theta}-u^*\|_{L_\rho^2(\Omega)}/\|u^*\|_{L_\rho^2(\Omega)}$ is $0.024$. The evolution of the
precision of the approximate solution is demonstrated in Figure~\ref{fig:testper_train}, which shows a
rapid convergence of the numerical solution. In order to compare the numerical solution with the ground truth, we plot the two functions along a randomly chosen axis in Figure~\ref{fig:randaxis}, from which we can see that the NN solution shows good accordance with the true solution. 

\subsection{Eigenvalue problem - Schr\"odinger operator}
In this section, we consider the eigenvalue problem associated with the Schr\"odinger operator
\begin{equation}
L = -\Delta + V,
\end{equation}
and the periodic boundary condition \eqref{eq:pbc}. Here $V$ is a potential function
\begin{equation}
V(x)=4\pi^2\sum_{i=1}^{d} c_{i} \cos \left(2\pi x_{i}\right),
\end{equation}
where $c_i\in[0, 0.2]$ for $1\leq i\leq d$. We adopt the same parameters as in
\cite{han2020solving}. The reference solution $u^*$ and the corresponding eigenvalue $\lambda^*$ are obtained by the spectral method described in
\cite{han2020solving}, and we measure the error of the numerical solutions by $E_0 = \|u_{\theta}-u^*\|_{L_\rho^2(\Omega)}/\|u^*\|_{L_\rho^2(\Omega)}$ and the error of the eigenvalue by $E_1 = |\lambda-\lambda^*|/|\lambda^*|$, where $\lambda$ is obtained via:
\[\frac{1}{\|u^*\|_{L_\rho^2(\Omega)}^2}\left(\frac{2}{\delta}\int_{\mathrlap{\Omega}} u_\theta(x)(u_\theta(x)-\E(u_\theta(x + W_\delta)))\mathrm{d} x + \int_{\mathrlap{\Omega}} V(x) |u_\theta(x)|^2 \mathrm{d} x\right),\]
if scheme I is used, and via
\[\frac{1}{\|u^*\|_{L_\rho^2(\Omega)}^2}\left(\frac{1}{\delta}\int_{\mathrlap{\Omega}} |u_\theta(x)-\E(u_\theta(x + W_\delta))|^2\mathrm{d} x + \int_{\mathrlap{\Omega}} V(x) |u_\theta(x)|^2 \mathrm{d} x\right).\]

Numerical tests are performed when $d=5, 10$. In the numerical results, $E_0$ is estimated on a test set of size $1\times 10^4$, and $E_1$ is estimated by averaging $10$ estimations on test sets of size $1\times 10^5$. The level $m$ of trigonometric bases used is set as $5$.

\subsubsection{Scheme I}
\begin{figure}[ht]
  \centering
  \subfigure[The error of the numerical solution in the process of training. \label{fig:err_5d1}]{
    \includegraphics[width=0.46\textwidth]{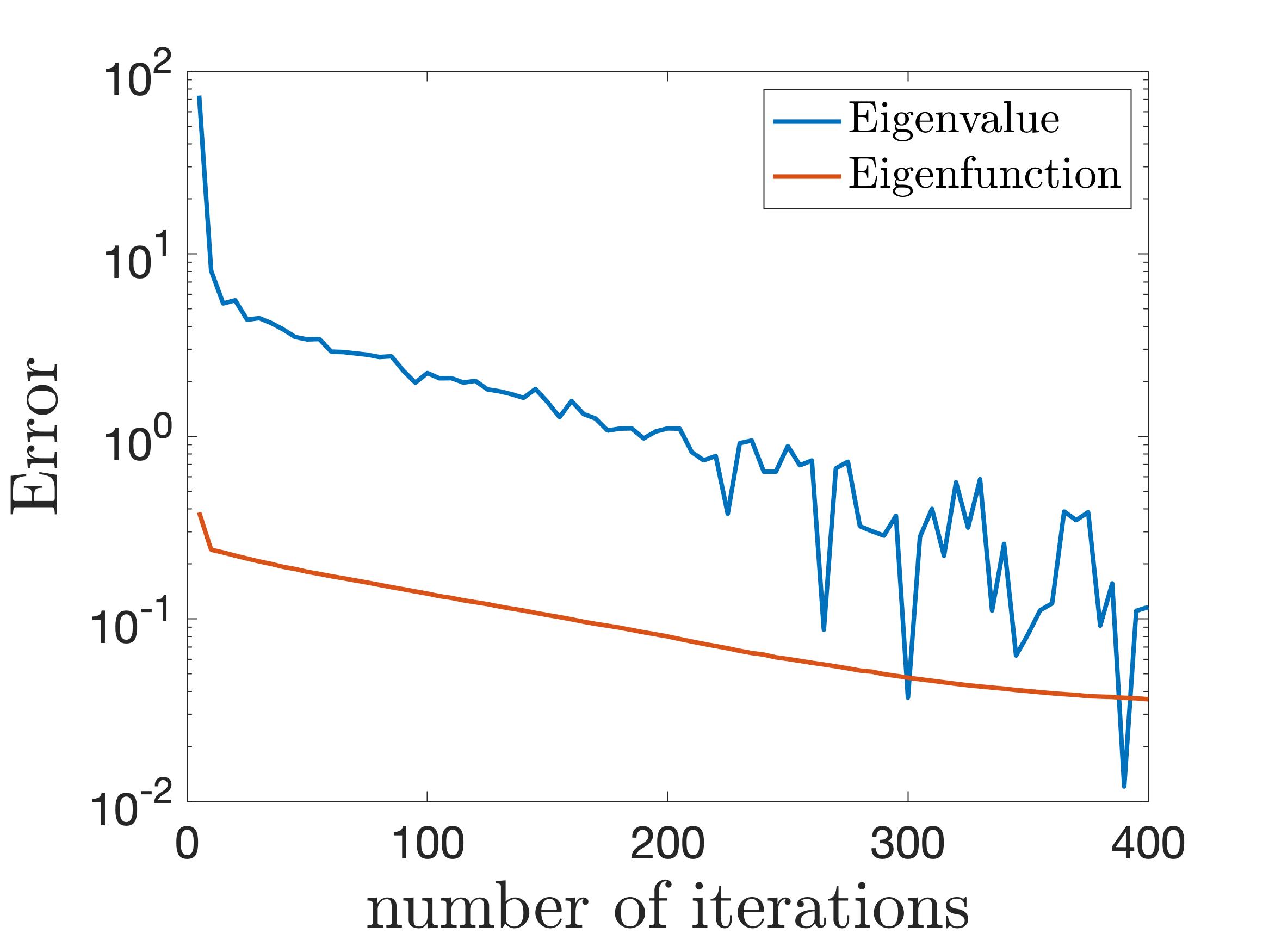} }
  \hspace{0.7em} \subfigure[The residue $\|u_\theta\|^2-1$ of the constraint in the process of training.\label{fig:norm5d1}]{
    \includegraphics[width=0.46\textwidth]{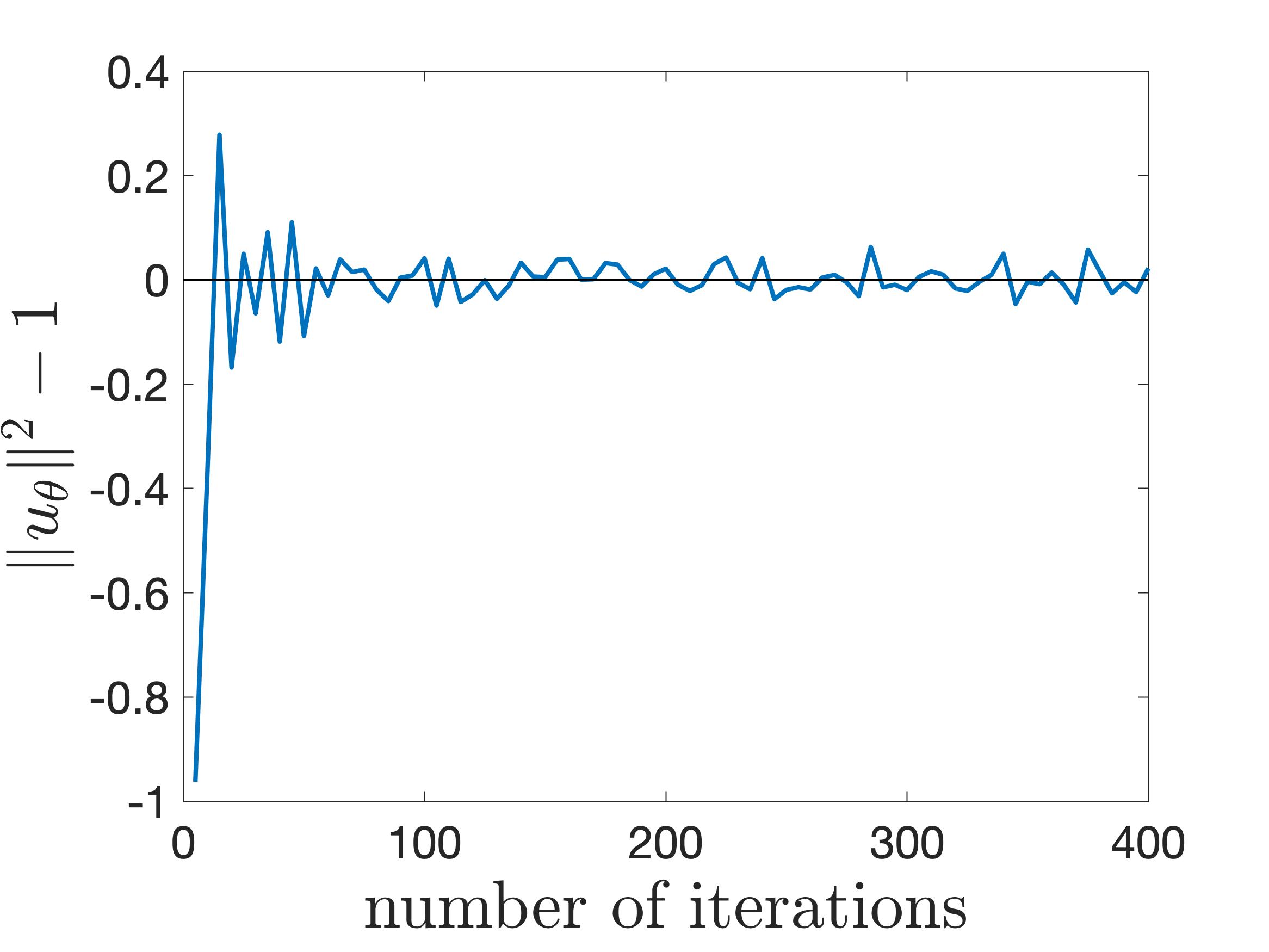} }
  \caption{The convergence of the error of the eigenfunction and eigenvalue and the residue of the constraint $\|u_\theta\|^2-1$ in the training process of the $5$-dimensional problem with scheme I. Blue curve in (a): Convergence of the error of the approximate eigenvalue; Orange curve in (a): Convergence of the error of the NN eigenfunction. }
  \label{fig:convergence5d1}
\end{figure}

When adopting the semigroup formulation \eqref{eq:eigenprimaldual1} to solve the $5$-dimensional case, $1.2\times 10^7$ samples are
used to train the neural network with width $300$, with $1.0\times 10^5$ samples used in each iteration. Throughout the training process, the learning rates $\tilde{\eta}_t=\eta_t=0.0003$. The hyper-parameters $g_{\text{default}}$, $c$ and $\delta$ are set to be $1$, $70$ and $0.0001$, respectively. After $400$ iterations of training, the final errors of the eigenfunction and eigenvalue are $E_0=0.036$ and $E_1 = 0.11$, respectively. The decay of errors is presented in Figure~\ref{fig:err_5d1}, from which we can see that the numerical solutions converge to the reference solution. 

\begin{figure}[ht]
  \centering \subfigure[The error of the numerical solution in the process of training. \label{fig:err_10d1}]{
    \includegraphics[width=0.46\textwidth]{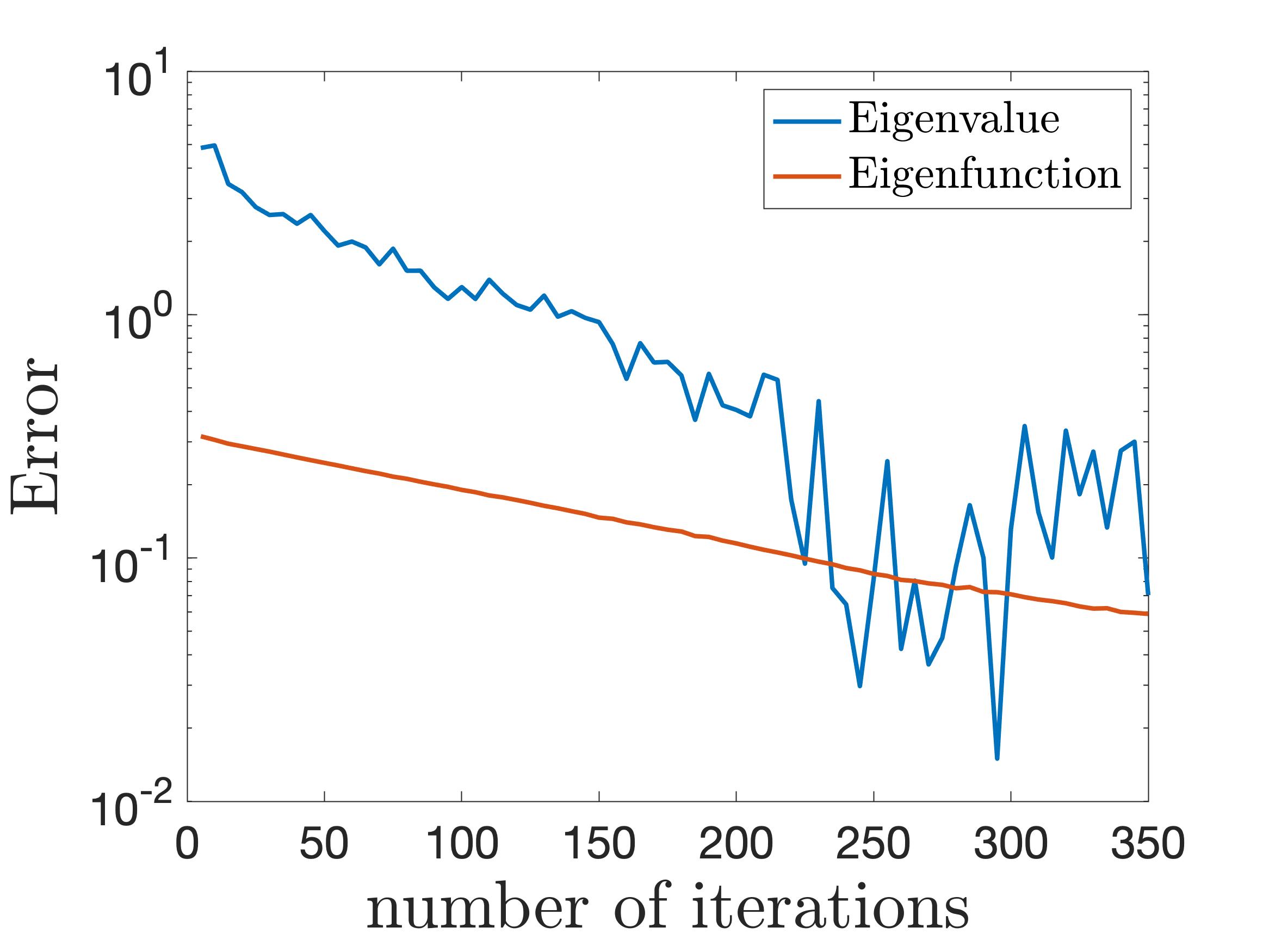} }
  \hspace{0.7em} \subfigure[The residue $\|u_\theta\|^2-1$ of the constraint in the process of training.\label{fig:norm10d1}]{
    \includegraphics[width=0.46\textwidth]{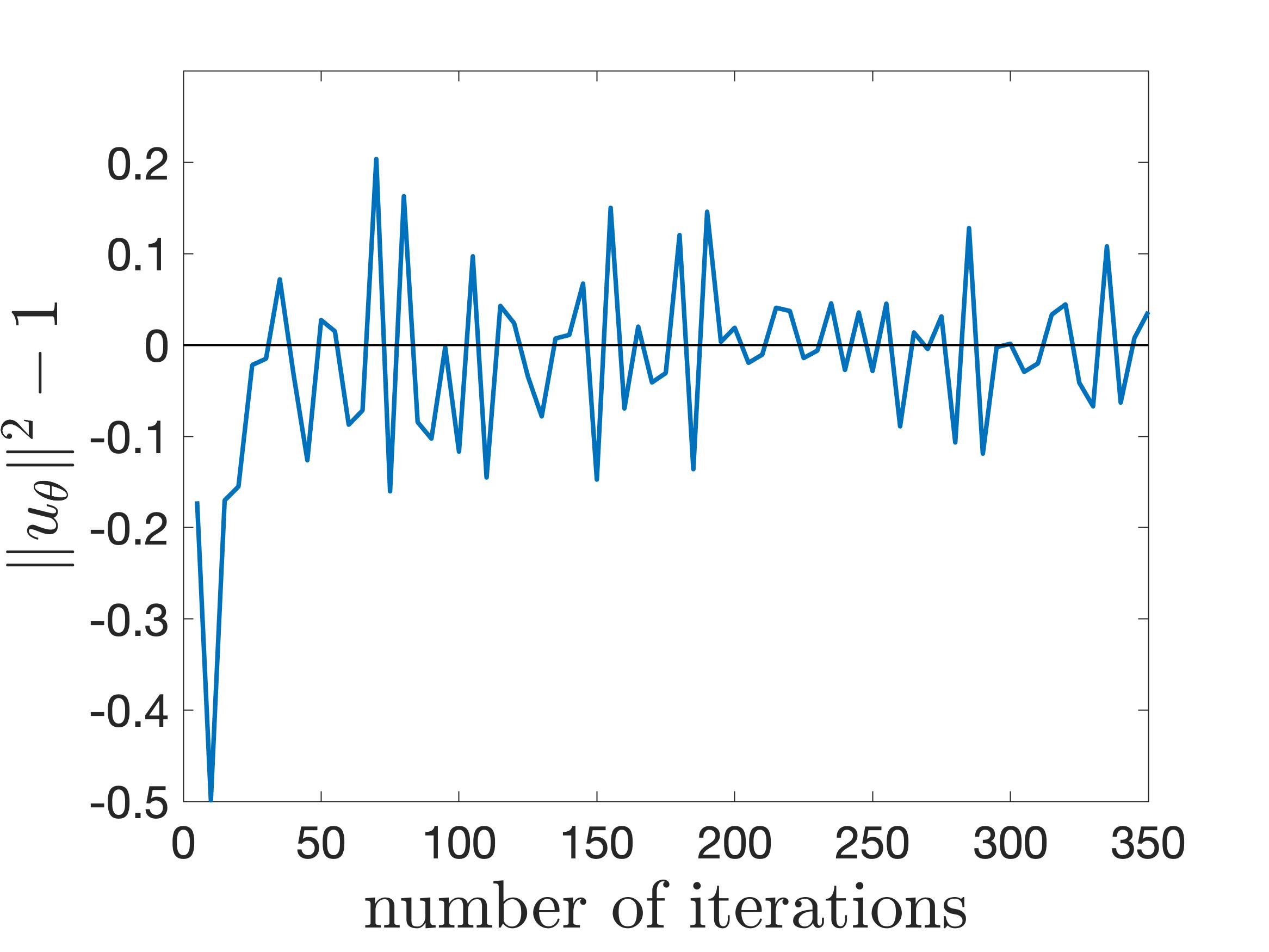} }
  \caption{The convergence of the error of the eigenfunction and eigenvalue and the residue of the constraint $\|u_\theta\|^2-1$ in the training process of the $10$-dimensional problem with scheme I. Blue curve in (a): Convergence of the error of the approximate eigenvalue; Orange curve in (a): Convergence of the error of the NN eigenfunction. } 
  \label{fig:convergence10d1}
\end{figure}

For the $10$-dimensional problem, we also use $1.2\times 10^7$ samples. A neural network with width $600$ is trained learning rates $\tilde{\eta}_t=\eta_t=0.0005$. In each iteration, $7.0\times10^4$ samples are randomly chosen from the training set. We set the hyper-parameters $g_{\text{default}}$, $c$ and $\delta$ as $4$, $40$ and $0.0001$, respectively. After $400$ iterations, the final errors of the eigenfunction and eigenvalue are $E_0=0.058$ and $E_1 = 0.07$, respectively. The convergence of the numerical solutions to the reference solution is shown in Figure~\ref{fig:err_10d1}. In both $5$d and $10$d problems, the constraint
$\|u_\theta\| = 1$ is well-enforced at the end of the training process, as shown in
Figure~\ref{fig:norm5d1} and Figure~\ref{fig:norm10d1}



\begin{figure}[ht]
  \centering \subfigure[Comparison of the 5d result. \label{fig:compare_5d1}]{
    \includegraphics[width=0.46\textwidth]{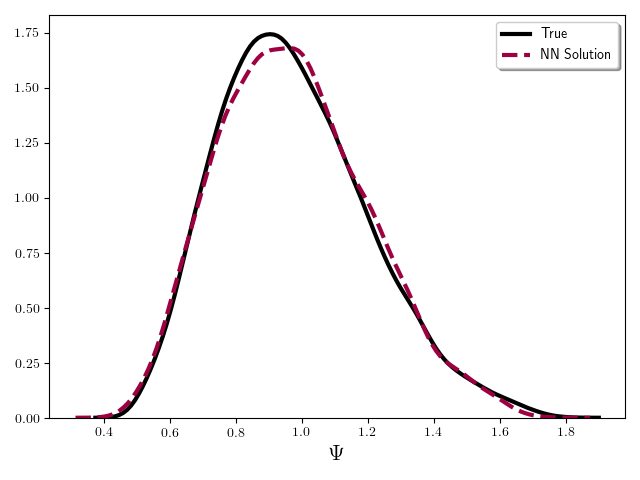} }
  \hspace{0.7em} \subfigure[Comparison of the 10d result.
  \label{fig:compare_10d1}]{
    \includegraphics[width=0.46\textwidth]{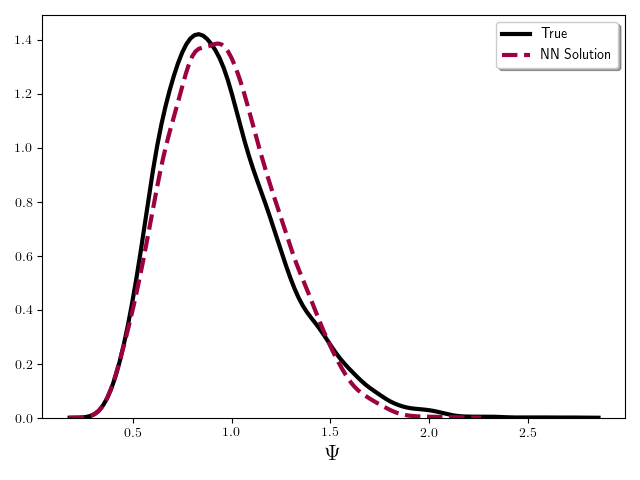} }
  \caption{Comparisons between the NN represented solutions and the ground truths using semigroup scheme I.
  (a): The 5-dimensional case.  (b): The 10-dimensional case. }
  \label{fig:compare1}
\end{figure}

Since it is difficult to visualize functions in high dimensions, we compare the probability density
function of $u_\theta(Z)$ and $u(Z)$, where $Z$ is a uniform random variable on $\Omega$. The
probability density function is obtained by performing kernel density estimation on a sample set of
size $10000$. The comparison of the probability density functions are given in
Figure~\ref{fig:compare1}, from which we can conclude that the numerical solutions obtained are in
good accordance with the ground truths.

\subsubsection{Scheme II}
When adopting the semigroup formulation \eqref{eq:eigenprimaldual2}, $4.0\times 10^6$ samples are used to train the neural network with width $300$ for the $5$-dimensional case and $600$ for the $10$-dimensional case, and $2000$ iterations of training is implemented. In each iteration, $1.0\times10^4$ samples are chosen randomly from the training set. In
this problem, the hyper-parameter $\delta$ and $c$ are set to be $0.001$ and $10$
respectively, and $g_{\text{default}}$ is set to be $4$ for the $5$-dimensional case and $1$ for the $10$-dimensional case. The learning rate for the dual variable
$\tilde{\eta}_t$ is set as $0.1$ while the learning rate for the primal variable $\eta_t$ is set as
$0.0008$ for the first half of the training process and $0.0003$ for the second half of the training
process. \\

\begin{figure}[ht]
  \centering \subfigure[The error of the numerical solution in the process of training. \label{fig:err_5d2}]{
    \includegraphics[width=0.46\textwidth]{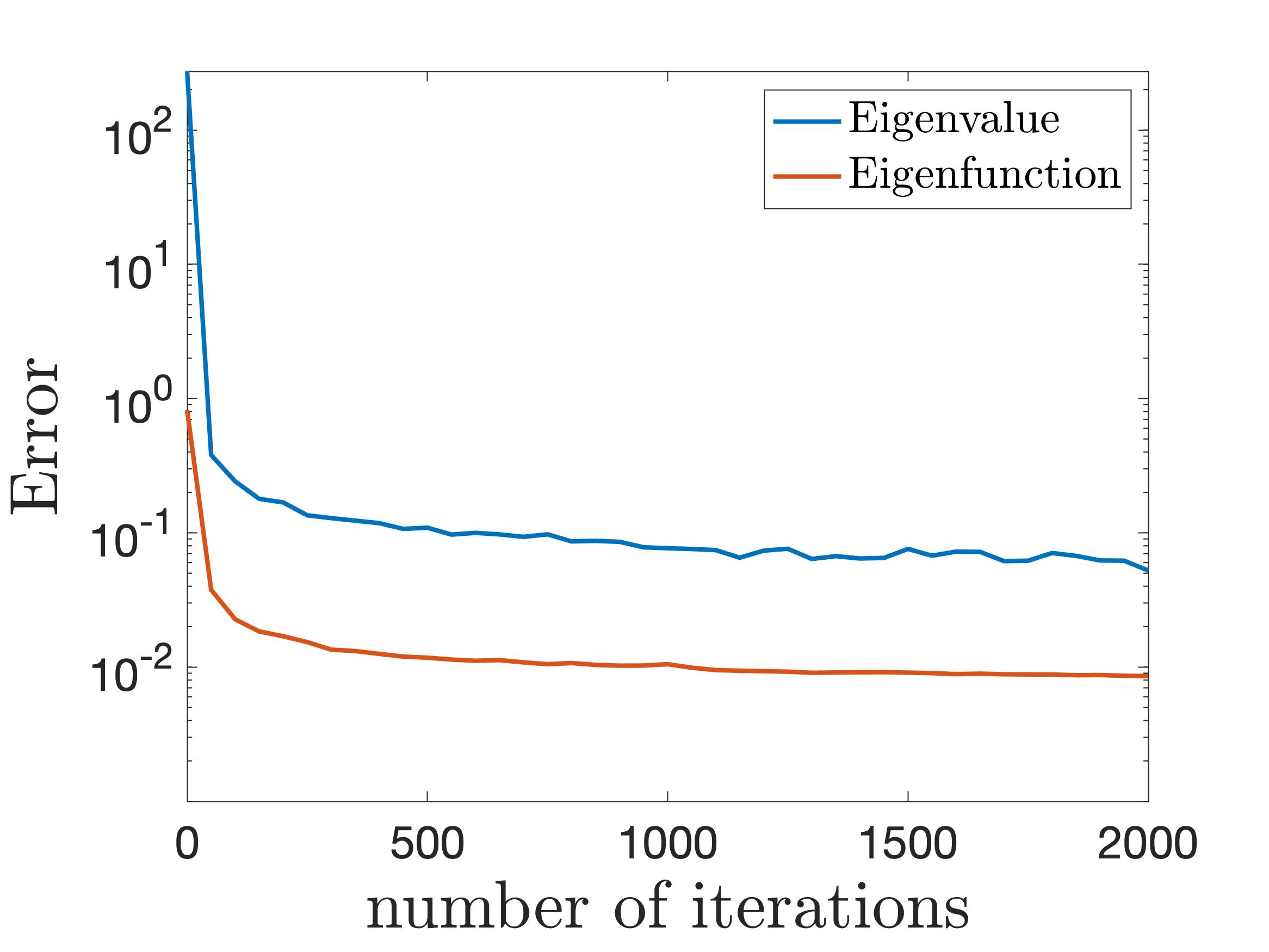} }
  \hspace{0.7em} \subfigure[The residue $\|u_\theta\|^2-1$ of the constraint in the process of training.\label{fig:norm5d2}]{
    \includegraphics[width=0.46\textwidth]{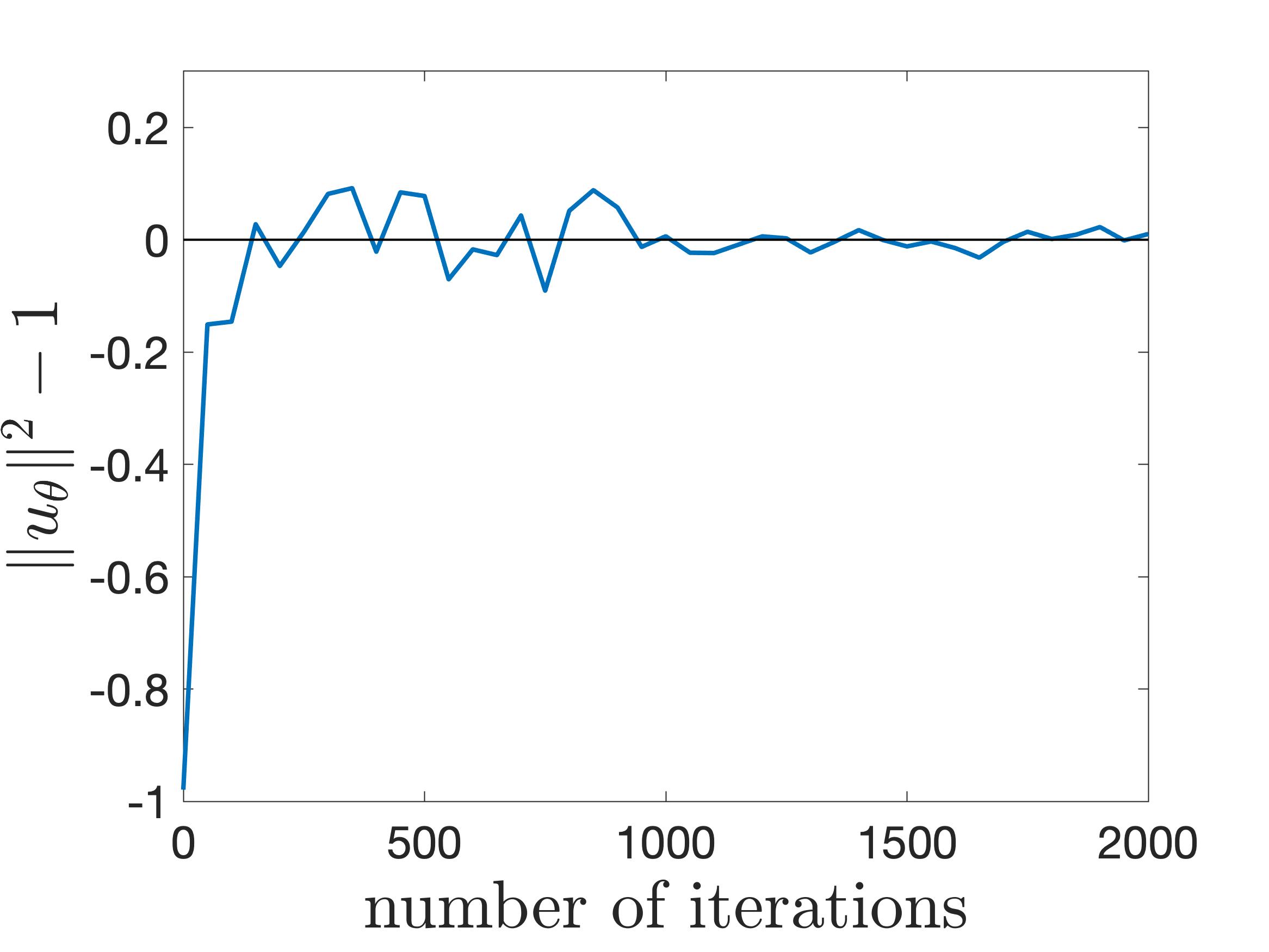} }
  \caption{The convergence of the error of the eigenfunction and eigenvalue and the residue of the constraint $\|u_\theta\|^2-1$ in the training process of the $5$-dimensional problem with scheme II. Blue curve in (a): Convergence of the error of the approximate eigenvalue; Orange curve in (a): Convergence of the error of the NN eigenfunction. } 
  \label{fig:convergence5d2}
\end{figure}

\begin{figure}[ht]
  \centering \subfigure[The error of the numerical solution in the process of training. \label{fig:err_10d2}]{
    \includegraphics[width=0.46\textwidth]{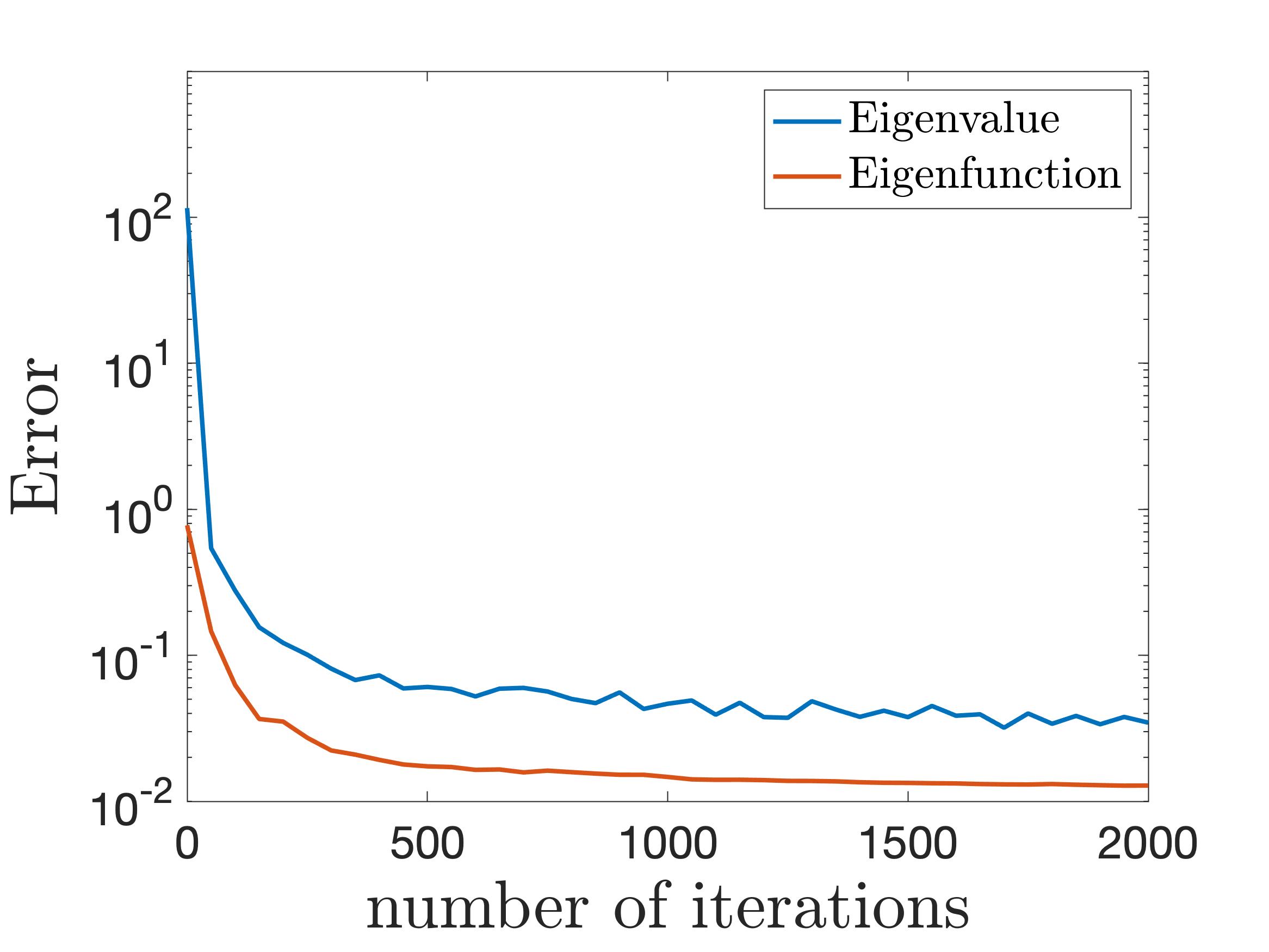} }
  \hspace{0.7em} \subfigure[The residue $\|u_\theta\|^2-1$ of the constraint in the process of training.\label{fig:norm10d2}]{
    \includegraphics[width=0.46\textwidth]{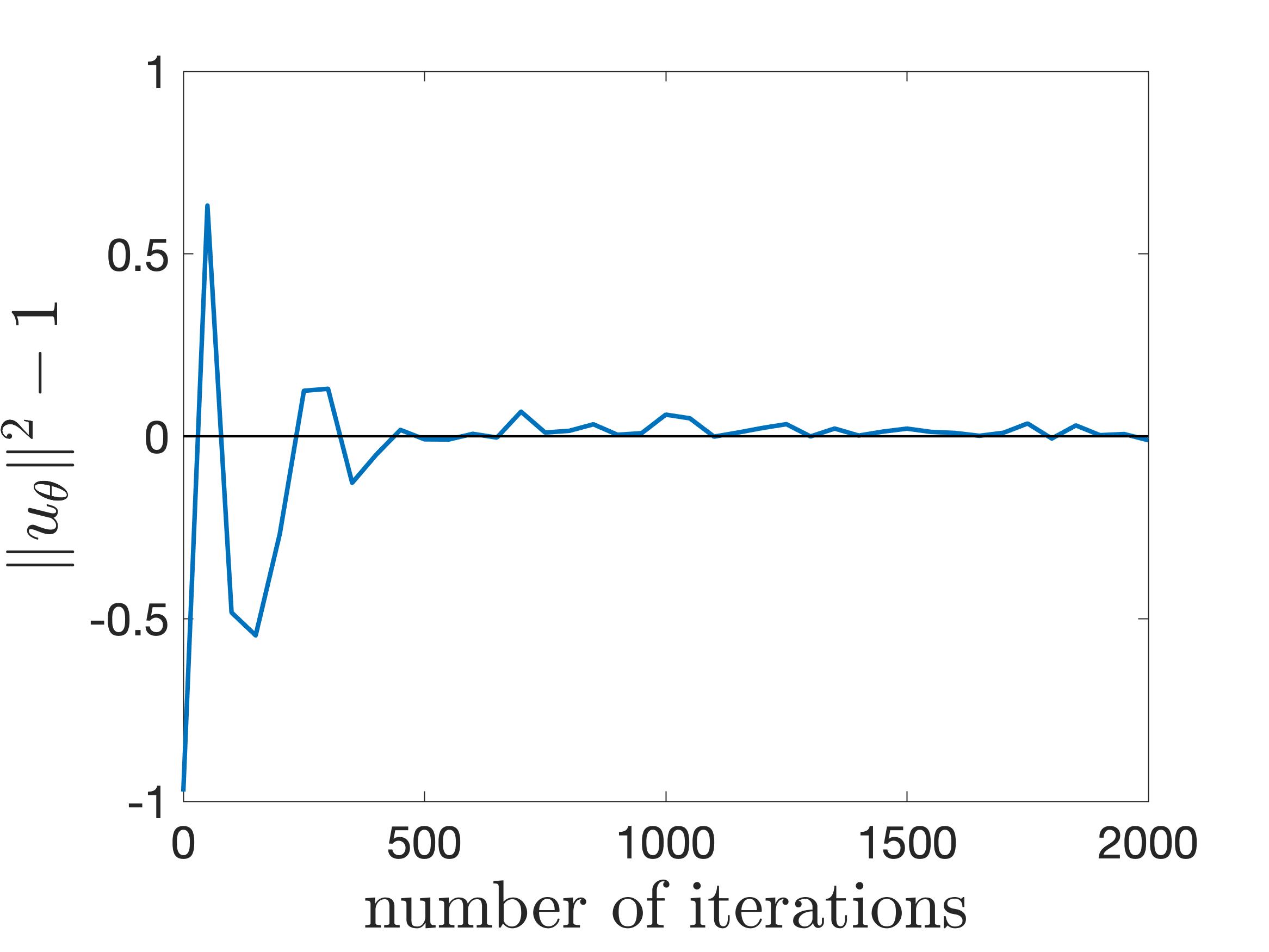} }
  \caption{The convergence of the error of the eigenfunction and eigenvalue and the residue of the constraint $\|u_\theta\|^2-1$ in the training process of the $10$-dimensional problem with scheme II. Blue curve in (a): Convergence of the error of the approximate eigenvalue; Orange curve in (a): Convergence of the error of the NN eigenfunction. }
  \label{fig:convergence10d2}
\end{figure}

\begin{figure}[ht]
  \centering \subfigure[Comparison of the 5d result. \label{fig:compare_5d2}]{
    \includegraphics[width=0.46\textwidth]{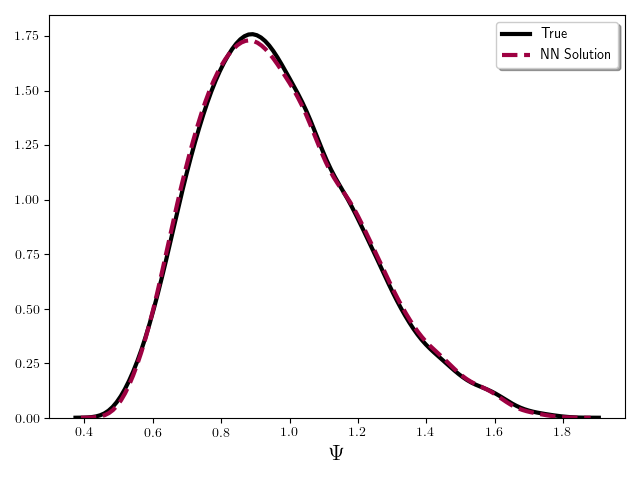} }
  \hspace{0.7em} \subfigure[Comparison of the 10d result.
  \label{fig:compare_10d2}]{
    \includegraphics[width=0.46\textwidth]{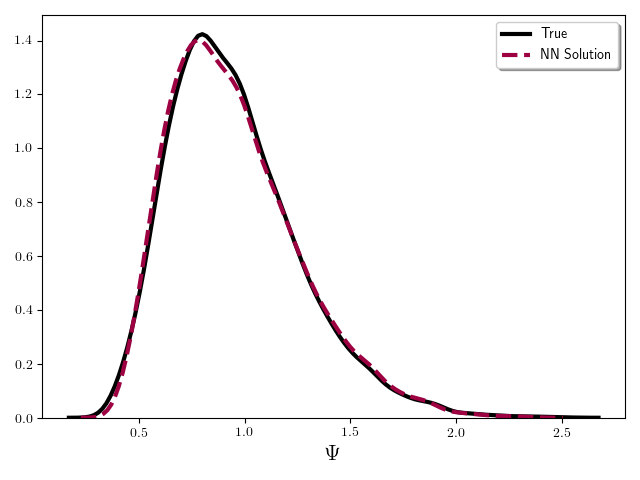} }
  \caption{Comparisons between the NN represented solutions and the ground truths using Scheme II.
  (a): The 5-dimensional case.  (b): The 10-dimensional case. }
  \label{fig:compare2}
\end{figure}

For the $5$-dimensional case and the $10$-dimensional case, the final errors of the eigenfunction are $0.0086$ and
$0.013$, respectively, and the final errors of the eigenvalue are $0.052$ and $0.035$, respectively. The training process are depicted in Figure~\ref{fig:convergence5d2} and
Figure~\ref{fig:convergence10d2}, respectively. As illustrated in Figure~\ref{fig:err_5d2} and
Figure~\ref{fig:err_10d2}, the numerical solutions of the eigenfunction and eigenvalue converge to the corresponding reference
solutions. Compared with the training process using Scheme I, the final errors are much
lower and the size of the training set is much smaller, although it takes more iterations to reach the final precision. This can also be verified by
comparing Figure~\ref{fig:compare1} and Figure~\ref{fig:compare2}, since the estimated probability
density functions for the numerical solutions are closer to the probability density functions for
the reference solutions in Figure~\ref{fig:compare2}. We check that the constraint $\|u_\theta\| =
1$ is well-enforced at the end of the training process in Figure~\ref{fig:norm5d2} and
Figure~\ref{fig:norm10d2}.

\textbf{Comparison with the method used in \cite{han2020solving}}
In this part, we briefly compare our numerical results obtained by scheme II with those in \cite{han2020solving}. The linear Schr\"odinger problem in their paper is only different with the problem here by a $2\pi$ factor. Specifically, they consider the operator
\[
L = -\Delta + V,
\]
and the periodic boundary condition \eqref{eq:pbc}, with the domain replaced by $[0, 2\pi]^d$, and $V$ replaced by
\[
V(x)=\sum_{i=1}^{d} c_{i} \cos \left(x_{i}\right),
\]
where the coefficients $\{c_i\}_{i=1}^d$ are the same with the coefficients used here. The $L^2$ error reported in \cite{han2020solving} is calculated after normalizing the solutions with $\int_\Omega u^2(x)\mathrm{d}x = |\Omega|$, which is also satisfied by the solutions here since the volume of $[0, 1]^d$ is $1$ and we assume $\|u^*\|_{L^2(\rho)}=1$. Since the $L^2$ error reported in their paper is invariant under scaling, it is reasonable to compare the numerical results obtained there with the results in this paper. In the following, we divide the errors of the eigenvalue reported in their paper by the true eigenvalues in order to be consistent with the measurement used here. 

We test the method in \cite{han2020solving} using the source code provided by the authors of \cite{han2020solving} and the hyperparameters in \cite{han2020solving} on the same machine as we used to test our methods, which is a machine with $4$ N$1$ virtual CPUs on the Google Cloud platform with altogether $26$ GB memory and a Tesla K80 GPU. In both the $5$-dimensional and $10$-dimensional cases, the model is trained for $80000$ iterations with $1024$ samples used in each iteration. For the $5$-dimensional case, the final errors for the eigenvalue and eigenfunction are $0.016$ and $0.0098$, respectively. For the $10$-dimensional case, the final errors for the eigenvalue and eigenfunction are $0.017$ and $0.012$, respectively. It can be seen that the precision of the eigenfunction is comparable with the precision of eigenfunction in scheme II, while the errors of the eigenvalue are smaller but of the same magnitude as our results. 

On the other hand, due to the simplicity of our method, our method enjoys a shorter computation time and require fewer samples. For the method in \cite{han2020solving}, the computation time used in the $5$-dimensional and $10$-dimensional problems are $1.5\times10^4$ seconds ($4.2$ hours) and $4.3\times10^4$ seconds ($12$ hours), respectively, and altogether $8.0\times 10^7$ samples are used, while for scheme II proposed in this paper the computation time used are $3.0\times 10^2$ seconds ($5.0$ minutes) and $5.7\times10^2$ seconds ($9.5$ minutes), and $4.0\times 10^6$ samples are used. If we consider the numerical solutions obtained by the first $11000$ iterations of the method in \cite{han2020solving}, the errors of the eigenvalue and eigenfunction are $0.087$ and $0.035$ for the $5$d case and $0.097$ and $0.059$ for the $10$d case, which is larger than those in our method, but the time used is still $2.1\times 10^3$ seconds ($35$ minutes) and $5.9\times10^3$ seconds ($98$ minutes) for the $5$d and $10$d case, respectively, and the number of samples used is still $1.1\times10^7$, which is much larger than our method. Finally, we mention that other problems such as the problem of finding the second eigenpair is addressed in \cite{han2020solving}, which is not covered in our paper. 
\section{Conclusion}

In this paper, we present a semigroup method solving high dimensional PDE problems and eigenvalue problems effectively. We have shown numerically the efficiency of the proposed method in problems that have non-zero right-hand-side term with Dirichlet boundary conditions and periodic boundary conditions. In comparison with popular deep PDE solvers such as the Deep Ritz method \cite{weinan2018deep} and PINN method \cite{raissi2019physics}, where penalty functions is used to enforce the Dirichlet boundary condition and thus changes the solution, the proposed semigroup method addresses the Dirichlet boundary condition without penalty functions, and even when penalty functions are used, the true solution remains the same. Two semigroup schemes are proposed for the eigenvalue problems. With a scalar Lagrange multiplier, these schemes are able to handle the constraint in the eigenvalue problem and obtain accurate solutions. In comparison with established solvers such as the BSDE method \cite{han2020solving}, the proposed method uses much less computation time to achieve the same precision for certain problem such as the linear Schr\"odinger problem. 

The numerical schemes adopted here are generally first-order schemes. For future work, higher order schemes can be applied to improve the precision. Moreover, importance sampling techniques can be integrated with the proposed method to facilitate the generation of a training set with high quality. For the PDE problems, we have extended the method in \cite{li2020solving} to the case with non-zero right-hand-side term and more general boundary conditions. We point out that it is possible to further generalize the semigroup method to other types of elliptic PDEs. One possible way to proceed is to replace the second-order derivatives using the semigroup operator in, for example, the PINN method (\cite{raissi2019physics}). In this way, it is possible to avoid the calculation of mixed third-order derivatives needed there and to integrate the boundary conditions naturally through the semigroup operator instead of enforcing it solely through a penalty function. We can possibly accelerate the training of the model and alleviate the difficulty in tuning the penalty coefficient by this replacement.

\bibliographystyle{abbrv}
\newpage
\bibliography{ref}


\end{document}